\documentclass[12pt]{amsart}
\usepackage{amsmath}
\usepackage{amsthm}
\usepackage{graphicx}
\usepackage{amssymb}
\usepackage{amscd}
\usepackage[abs]{overpic}        
\newcommand{\Mob}{\mathrm{M\ddot{o}b}}

\newcommand{\Isom}{\mathrm{Isom}}
\newcommand{\Stab}{\mathrm{Stab}}
\newcommand{\Fix}{\mathrm{Fix}}
\newcommand{\sm}{\setminus}  
\newcommand{\fty}{\infty}

\newcommand{\bb}{\mathbb}  
  
\newcommand{\mcal}{\mathcal}  
\newcommand{\psl}{{\mathrm{PSL}}_2({\bb C})}

\newcommand{\id}{{\mathrm{id}}}
\newcommand{\bd}{\partial}
  
\newcommand{\im}{\mathrm{Im}\,}     
\newcommand{\re}{\mathrm{Re}\,} 
\newcommand{\ov}{\overline}
  
\newcommand{\bs}[1]{\boldsymbol{#1}}

\newcommand{\R}{\mathbb{R}}
\newcommand{\C}{\mathbb{C}}
\newcommand{\Z}{\mathbb{Z}}
\newcommand{\N}{\mathbb{N}}
\newcommand{\M}{\mcal M}

\newcommand{\bsv}{{\bs v}}
\newcommand{\bsu}{{\bs u}}
\newcommand{\bse}{{\bs e}}
\newcommand{\bsx}{{\bs x}}
\newcommand{\bsa}{{\bs a}}

\newcommand{\bsp}{{\bs p}}
\newcommand{\wh}{\widehat}

\newcommand{\hR}{\hat \R} 
\newcommand{\two}{\Mob(\hR^2)}
\newcommand{\three}{\Mob(\hR^3)}
\newcommand{\la}{\langle}  
\newcommand{\ra}{\rangle} 
\newcommand{\B}{{\bf B}} 
\newcommand{\HH}{{\bf H}} 
\newcommand{\D}{{\bf D}} 
\newcommand{\wa}{\mathring} 

\newcommand{\Si}{\Sigma} 
 
\newcommand{\ac}{\acute} 
 
\newcommand{\hC}{\hat \C} 
\newcommand{\ch}{\check} 
 
\newcommand{\hP}{\hat P} 
\newcommand{\CC}{\mathcal C} 
\newcommand{\MM}{\M_{1,1}} 
\newcommand{\MN}{\M_{0,4}} 
\newcommand{\hMM}{\wh \M_{1,1}} 
\newcommand{\radi}{\mathrm{radi}} 
\newcommand{\DD}{{\mathcal D}} 
\newcommand{\NN}{{\mathcal N}} 
\theoremstyle{plain}
  \newtheorem{thm}{Theorem}[section]
  \newtheorem{cor}[thm]{Corollary}
  \newtheorem{lem}[thm]{Lemma}

\theoremstyle{definition}
  \newtheorem{defn}[thm]{Definition}
  
\theoremstyle{remark}
  \newtheorem*{rem}{Remark}
  
  \newtheorem*{ack}{Acknowledgements} 

\numberwithin{equation}{section}
\title[An extension of the Maskit slice]
{An extension of the Maskit slice for $4$-dimensional Kleinian groups}

\author{Yoshiaki Araki}
\address{Synclore Corporation, Hakuyo Building, 3-10 Nibancho Chiyoda-ku, 
Tokyo 102-0084, Japan}
\email{yoshiaki.araki@synclore.com}

\author{Kentaro Ito}
\address{ Graduate School of Mathematics, 
Nagoya University, Nagoya 464-8602, Japan}
\email{itoken@math.nagoya-u.ac.jp}

\date{October 8, 2008}

\begin{document}
\begin{abstract}
Let $\Gamma$ be a $3$-dimensional Kleinian punctured torus group 
with accidental parabolic transformations. 
The deformation space of $\Gamma$ in the group 
of M\"{o}bius transformations on the $2$-sphere 
is well-known as the Maskit slice $\MM$ of punctured torus groups.  
In this paper, we study deformations $\Gamma'$
of $\Gamma$ in the group of M\"{o}bius transformations on the $3$-sphere 
such that $\Gamma'$ does not contain screw parabolic transformations. 
We will show that the space of the deformations 
is realized as a domain of $3$-space $\mathbb{R}^3$, 
which contains the Maskit slice $\MM$ as a slice through a plane. 
Furthermore, we will show that 
the space also contains the Maskit slice $\MN$ 
of fourth-punctured sphere groups as a slice through another plane. 
Some of another slices of the space will be also studied. 
\end{abstract}

\maketitle

\section{Introduction}

Let $\Gamma$ be a $3$-dimensional Kleinian once-punctured (or simply punctured) 
torus group with accidental parabolic transformations, 
which acts on the Riemann sphere $\hC$ as M\"{o}bius transformations. 
In this paper, we regard $\Gamma$ as a $4$-dimensional Kleinian group 
and study deformations $\Gamma'$ of $\Gamma$ in the group of 
M\"{o}bius transformations on the $3$-sphere $\hR^3=\R^3 \cup \{\fty\}$; 
a Kleinian group $\Gamma'$ is called a {\it deformation} of $\Gamma$ if 
there is an isomorphism $\phi:\Gamma \to \Gamma'$ which takes a 
parabolic transformation to a parabolic transformation. 
Especially, we focus on all deformations $\{\phi:\Gamma \to \Gamma'\}$ of $\Gamma$ such that 
$\phi$ takes  a pure parabolic transformation to a pure parabolic transformation; 
the definition of a pure/screw parabolic transformation will be described below (see also 2.5). 

The $n$-dimensional sphere $\hR^n=\R^n \cup \{\fty\}$ 
is naturally identified with the ideal boundary of the 
$(n+1)$-dimensional hyperbolic space $\HH^{n+1}$. 
Therefore, the group $\Mob(\hR^n)$ of 
orientation preserving M\"{o}bius transformations of $\hR^n$ is 
identified with the group $\Isom(\HH^{n+1})$ of 
 orientation preserving isomorphisms of $\HH^{n+1}$. 
A discrete subgroup of $\Isom(\HH^{n+1})=\Mob(\hR^n)$ is said to be a 
$(n+1)$-dimensional {\it Kleinian group}. 

$3$-dimensional Kleinian groups are well studied 
in various contexts; for example, hyperbolic geometry, $3$-dimensional topology, 
complex analysis, etc. 
We refer the reader to Kapovich \cite{Ka2} and Marden \cite{Mar} for overviews of the theory of 
$3$-dimensional Kleinian groups. 
 On the other hand, the study of 
Kleinian groups in higher dimensions is more wild and is 
far from getting the whole picture. 
We refer the reader to Apanasov \cite{Ap} and 
Kapovich \cite{Ka1} for more information of this area. 
Among higher dimensional Kleinian groups, the study of 
$4$-dimensional Kleinian groups should be 
interesting in its own light 
because they have ``fractal" limit sets in the $3$-sphere $\hR^3$ into which our geometric insight work well. 
Some examples of $4$-dimensional Kleinian groups 
with computer graphics of their limit sets can be found in Ahara--Araki \cite{AA}. 

We now move into more detailed argument. 
In what follows, we regard $\two$ as a subgroup of $\three$ by embedding 
$\hR^2$ into $\hR^3$ by $(x,y) \mapsto (x,y,0)$.  
Recall that the group $\two$ 
is naturally identified with the group $\mathrm{Aut}(\hC)$ of 
conformal automorphisms of the Riemann sphere $\hC$, 
or the group of linear fractional transformations.  
Let $\Gamma \subset \two$ be a 
{\it Kleinian punctured torus group with accidental parabolic transformations}; 
that is, $\Gamma$ is a $3$-dimensional Kleinian group such that 
the quotient manifold $\HH^3/\Gamma$ is homeomorphic to 
the trivial interval bundle over a punctured torus, 
and the ideal conformal boundary of $\HH^3/\Gamma$
is a union of a punctured torus and a thrice-punctured sphere. 
Then $\Gamma$ is a rank-$2$ free group and a pair of generators $\alpha, \,\beta \in \two$ of $\Gamma$ 
can be chosen so that $\beta$ and the commutator $[\alpha,\beta]$ of the generators are parabolic. 
The following normalization of such a group $\Gamma$ up to conjugation is fundamental. 

\begin{lem}[cf. \cite{Kr}, \cite{KS}]
Let $\Gamma=\la \alpha, \beta \ra$ 
be a (not necessarily discrete) rank-$2$ free group in $\two$ 
such that $\beta$ and $[\alpha,\beta]$ are parabolic.
Then there is a complex number $\mu \in \C$ such that 
$\Gamma$ is conjugate in $\two$ to the group 
$G_\mu=\la A_\mu,B\ra$ 
generated by two M\"{o}bius transformations 
\begin{eqnarray}
A_\mu(\tau)=\frac{1}{\tau}+\mu \quad \text{and} \quad 
B(\tau)=\tau+2.  
\end{eqnarray}
\end{lem}

We define: 
\begin{eqnarray*}
\MM=\{\mu \in \C : G_\mu=\la A_\mu,B\ra \ 
\text{is a rank-$2$ free Kleinian group}\}. 
\end{eqnarray*}
This set is known as the {\it Maskit slice} of $3$-dimensional 
Kleinian punctured torus groups (see \cite{KS} for more information). 

In this paper we want to consider deformations of a $3$-dimensional Kleinian group 
$G_\mu=\la A_\mu, B\ra$ with $\mu \in \MM$ in the group $\three$. 
To explain our results, we need to recall a classification of elements of $\three$; 
see Section 2 for more details. 
Similar to the case of $\two$, elements of $\three$ are classified into three types:
 elliptic, parabolic and loxodromic transformations. 
Especially a transformation in $\three$ is {\it parabolic} if it has exactly one fixed point in $\hR^3$. 
Those parabolic elements are further classified into two types: 
a parabolic transformation $f \in \three$ is said to be {\it pure parabolic} 
if it is conjugate in $\three$ to a translation,  
and {\it screw parabolic} otherwise. 
A screw parabolic transformation is conjugate to a composition 
of a rotation and a translation along a common axis. 

The first goal of this paper is to show the following theorem which is 
an analogue of Lemma 1.1 for groups in $\three$:  

\begin{thm}(Theorem 4.5)
Let $\Gamma=\la \alpha, \beta \ra$ 
be a (not necessarily discrete) rank-$2$ free group in $\three$ 
such that $\beta$ and $[\alpha,\beta]$ are pure parabolic.
Then there is a point $\bsp=(p,q,r)$ in $\R^3$ such that 
$\Gamma$ is conjugate in $\three$ to the group 
$G_\bsp=\la A_\bsp,B\ra$ 
generated by 
\begin{eqnarray*}
&&A_\bsp(x,y,z)=\frac{(x,-y,z)}{x^2+y^2+z^2}+(p,q,r) \quad \text{and} \\
&&B(x,y,z)=(x,y,z)+(2,0,0). 
\end{eqnarray*}
\end{thm}

It is just a technical reason why we assume in the theorem above that 
$\beta$ and $[\alpha,\beta]$ are ``pure parabolic" instead of ``parabolic,"   
and it is more natural to consider the case where 
$\beta$ and $[\alpha, \beta]$ could be screw parabolic. 
In fact,  suppose that a homeomorphism $f:\hR^3 \to \hR^3$ 
conjugates a $3$-dimensional Kleinian group $\Gamma \subset  \two$ 
with  a parabolic transformation $\gamma \in \Gamma$ 
to a  $4$-dimensional Kleinian group $\Gamma'=f\Gamma f^{-1} \subset \three$. 
Then  the transformation $f \gamma f^{-1} \in \Gamma'$ 
could be screw parabolic in general (see 2.5).   
Unfortunately, we do not know the result of Theorem 1.2 
with ``pure parabolic" replaced by ``parabolic."   
We also remark that, recently, Y.~Kim and K.~Sakugawa independently announced that 
the level $2$ congruence subgroup of the modular group, 
the rank-$2$ free group $\Gamma=\la \alpha,\beta \ra$ in $\two$ 
such that $\alpha,\beta$ and $\alpha\beta$ are parabolic, 
has a continuous family of deformations $\{\Gamma'\}$ in $\three$ 
such that $\Gamma'$ are Kleinian groups containing screw parabolic transformations 
(see also Theorem 4.2 and its remark). 

We define: 
\begin{eqnarray*}
\hMM=\{\bsp=(p,q,r) \in \R^3 : G_\bsp=\la A_\bsp,B\ra \ 
\text{is a rank-$2$ free Kleinian group}\}. 
\end{eqnarray*}
The next purpose of this paper is to understand  
the shape of the set $\hMM$. 
We remark that if $\bsp=(p,q,0)$, the group $G_\bsp=\la A_\bsp,B \ra$ is 
nothing but the Poincar\'{e} extension of the group 
$G_\mu=\la A_\mu, B \ra$ with $\mu=p+iq$. 
Therefore $\hMM$ contains $\MM$ as a slice through the plane $r=0$ in the 
parameter space $\R^3=\{\bsp=(p,q,r)\}$. 
On the other hand, 
we will show in Section 5 the slice of $\hMM$ through 
the plane $q=0$ is the Maskit slice $\MN$ of $3$-dimensional 
fourth-punctured sphere groups which is defined by 
\begin{eqnarray*}
\MN=\{\mu \in \C : H_\mu=\la B, C,D_\mu \ra \ 
\text{is a rank-$3$ free Kleinian group}\}, 
\end{eqnarray*}
where $B(\tau)=\tau+2$, $C(\tau)=\tau/(2\tau+1)$ 
and $D_\mu(\tau)=C(\tau-\mu)+\mu$ (see Section 3).  

\begin{thm}(Theorems 5.4 and 5.7)
The Maskit slice $\MN$ appears as the slice of $\hMM$ through the plane $q=0$. 
Moreover, there is a constant $0<\varphi_0<\pi/2$ such that 
for every $\theta \in [-\varphi_0,  \varphi_0]$, 
$\MN$ also appears as 
the slice through the plane 
$q=0$ with angle $\theta$ rotated along the $p$-axis.  
\end{thm}

Although the boundary of the slices of $\hMM$ 
through the planes 
 $r=0$ and $q=0$ are so called ``fractal," 
we will also obtain in Section 5 
the following theorem (see also Figure 1).    

\begin{thm}(Corollary 5.11) 
The boundary of the slice of $\hMM$ through the plane $p=0$ is 
a union of countably many analytic arcs.    
\end{thm}

\begin{figure}[htbp]
  \begin{center}
      \includegraphics[keepaspectratio=true,height=70mm]{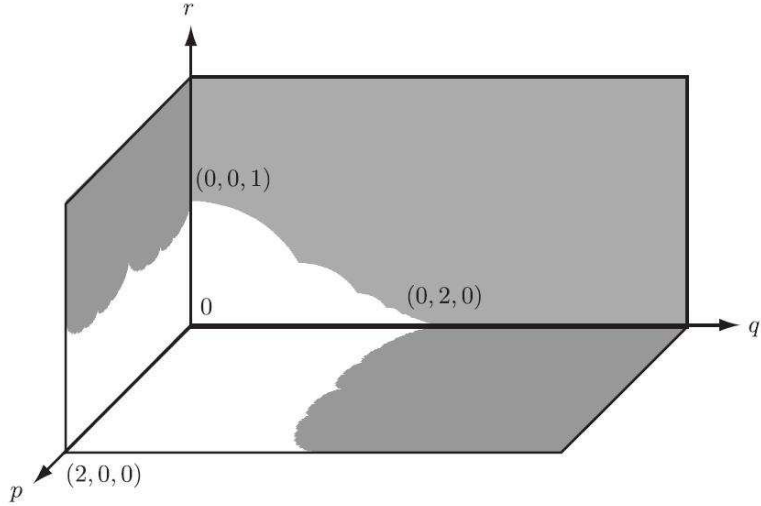}
  \end{center}
  \caption{A schematic figure of slices of $\hMM$ through the planes 
  $p=0$, $q=0$ and $r=0$.}
\end{figure}
One can find in Figure 2 some computer graphics 
of the limit sets of the groups $G_\bsp$ for parameters 
$\bsp \in \R^3$ which lie (or seem to be lie) in the set $\hMM$. 
These figures can be seen as $3$-dimensional extensions 
of the beautiful patterns in the complex plane included in the book \cite{MSW} 
by Mumford, Series and Wright. 
\begin{figure}[thbp]
  \begin{center}
    \includegraphics[keepaspectratio=true,height=190mm]{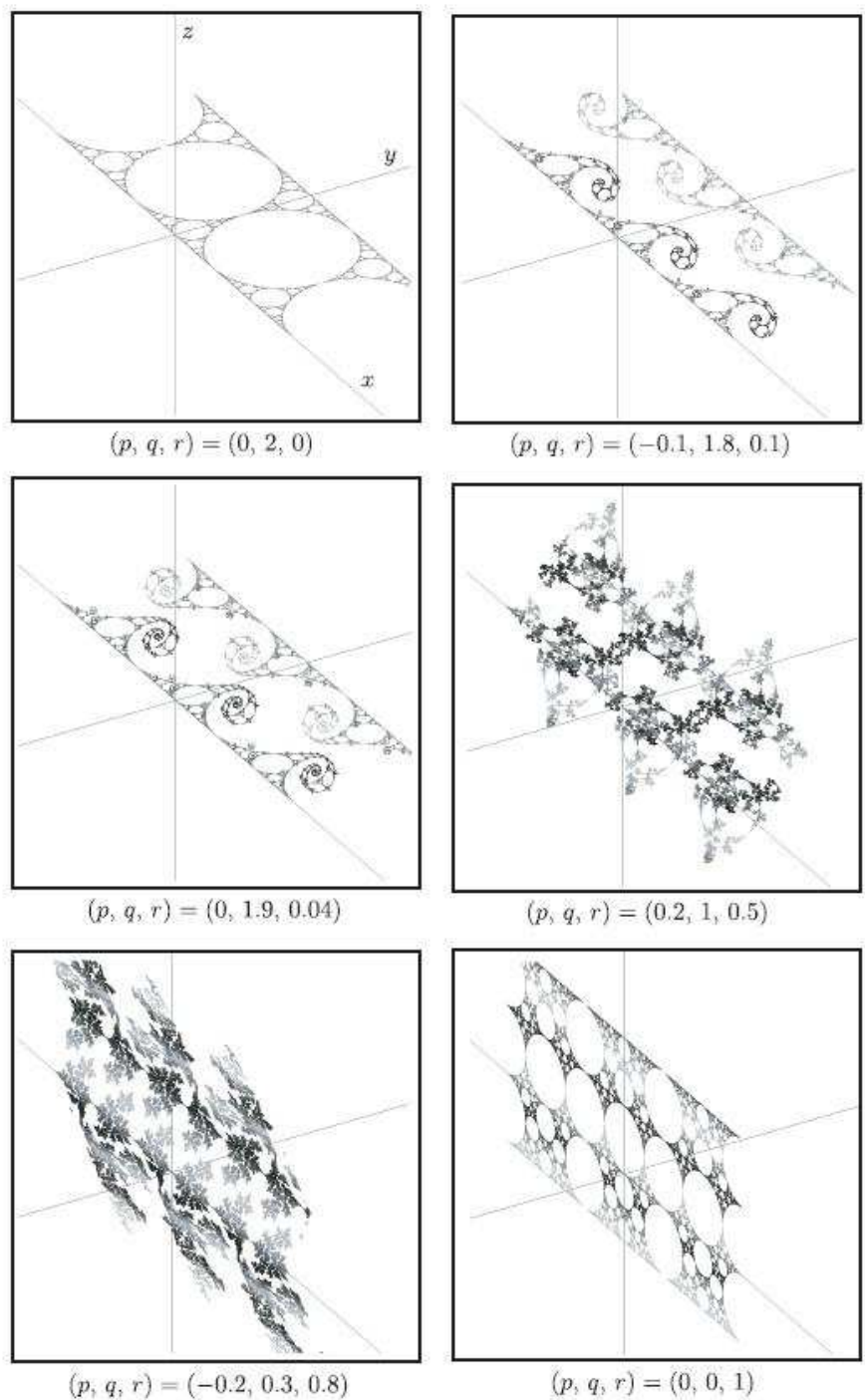}
  \end{center}
  \caption{Computer graphics of the limit sets  
$\Lambda(G_\bsp)$ of $G_\bsp$ restricted to the domain 
$\{(x,y,z) \in \R^3 : |x| \le 3\}$. The brightness of a point
$(x,y,z) \in \Lambda(G_\bsp)$ depends on its $z$-value. 
These are the views form the same point.}
\end{figure}
 
We remark that M\"{o}bius transformations in  $\three$ 
can be related to two-by-two matrices whose entries lie in the quaternion algebra,  
see Cao--Parker--Wang \cite{CPW}. 
However, in this paper, we take a geometric approach 
and make use of geometric techniques in the theory of Kleinian groups, 
instead of calculating in the quaternion algebra. 

This paper is organized as follows: 
In Section 2, we recall some basic facts of M\"{o}bius transformations.   
In Section 3,  we 
recall the definition of the 
Maskit slice $\MM$ (resp. $\MN$) of $3$-dimensional 
punctured torus groups (resp. fourth-punctured sphere groups).  
In section 4, we prove Theorem 1.2 on a normalization of $4$-dimensional 
punctured torus groups with accidental parabolics. 
In section 5, we study the shape of the space $\hMM$ and prove 
Theorems 1.3 and 1.4. 
As an application of Theorem 1.2,  
we give in Appendix an example of a family of $4$-dimensional Kleinian groups 
with $3$-generators, which contain 
punctured torus groups as $2$-generator subgroups. 

\begin{ack}
The authors would like to thank 
Shin Nayatani and Yasushi Yamashita for their interest and useful comments, 
Yohei Komori for his valuable information and discussions, 
and Keita Sakugawa for bringing to our attention 
to deformations of Kleinian groups with screw parabolic transformations.  
We also grateful for the referee's valuable comments and suggestions. 
\end{ack}

\section{Preliminaries}

In this section, we recall from Beardon \cite{Be} and  Matsumoto \cite{Mat} 
the basic facts of M\"{o}bius transformations on 
$n$-dimensional sphere $\hR^n=\R^n \cup \{\fty\}$. 
In the succeeding sections, we are mainly concerned with the cases of $n=2,\,3$.  

\subsection{Inversions.} 
For $\bsx=(x_1, \ldots, x_n) \in \R^n$, 
$$
|\bsx|=\sqrt{x_1^2+\cdots+x_n^2}
$$
denotes the Euclidean norm of $\bsx$. 
A $(n-1)$-dimensional {\it sphere} $\sigma$ in $\hR^n$ is either 
an Euclidean sphere $\{\bsx \in \R^n: |\bsx-\bsa|=r\}$ with 
$\bsa \in \R^n$, $r>0$,   
or a $(n-1)$-dimensional 
Euclidean plane $P$ plus $\{\fty\}$. 
The {\it inversion} $J_\sigma:\hR^n \to \hR^n$ in the sphere 
$\sigma$ is defined as follows: 
if $\sigma=\{\bsx \in \R^n: |\bsx-\bsa|=r\}$   
then   
\begin{eqnarray}
J_{\sigma}(\bsx)=r^2\frac{\bsx-\bsa}{|\bsx-\bsa|^2}+\bsa, 
\end{eqnarray}
and if 
$\sigma=P \cup \{\fty\}$ then $J_\sigma$ 
is the reflection in the plane $P$.  
Especially 
we denote the inversion in the unit sphere at the origin by 
$$
J(\bsx)= \frac{\bsx}{|\bsx|^2}. 
$$
Then the inversion $J_\sigma$ in $(2.1)$ is written as 
\begin{eqnarray}
J_\sigma(\bsx)=r^2J(\bsx-\bsa)+\bsa. 
\end{eqnarray}

\subsection{M\"{o}bius transformations.}  
A {\it M\"{o}bius transformation} on $\hR^n$ 
is a product of finite numbers of 
inversions in codimension one spheres in $\hR^n$. 
We denote by $\Mob(\hR^n)$ the group of orientation preserving  
M\"{o}bius transformations. 
It is known that the group $\Mob(\hR^n)$ 
is equal to the group of all orientation preserving 
conformal automorphism of $\hR^n$. 
Especially, the group $\two$ is naturally identified with the group 
$\mathrm{Aut}(\hat \C)$  of orientation preserving 
conformal automorphism of $\hat \C$ 
and with the Lie group $\psl$. 
The followings are typical examples of transformations $f:\hR^n \to \hR^n$ 
in $\Mob(\hR^n)$: 
\begin{itemize}
  \item A {\it translation}: $f(\bsx)=\bsx+\bsa$,   
  where $\bsa \in \R^n$. We denote by  $T_\bsa(\bsx)=\bsx+\bsa$. 
  \item A {\it magnification}:  $f(\bsx)=\lambda \bsx$,  where $\lambda>0$. 
  \item An {\it orthonormal transformation}: $f(\bsx)=P(\bsx)$,  where  
  $P \in SO(n)$. 
\end{itemize}
General M\"{o}bius transformations can be written as follows. 
\begin{lem}[Theorem 3.5.1 in \cite{Be}]
Let $f \in \Mob(\hR^n)$. 
\begin{enumerate}
  \item If $f(\fty)=\fty$ then $f(\bsx)=\lambda P(\bsx)+\bsu$, where 
$\lambda >0$, $P \in SO(n)$ and $\bsu \in \R^n$. 
  \item If $f(\fty) \ne \fty$ then 
$f(\bsx)=\lambda PJ (\bsx-\bsu)+\bsv$, where 
$\lambda >0$, $P \in O(n) \sm SO(n)$ and $\bsu, \bsv \in \R^n$. 
\end{enumerate}
\end{lem}

\subsection{Isometric spheres.}  
Let $f \in \Mob(\hR^n)$ be the same as in Lemma 2.1 (2). 
Note that $\bsu=f^{-1}(\fty)$ and $\bsv=f(\fty)$. 
The {\it isometric sphere} $I(f)$ of $f$ is defined by 
$$
I(f)=
\{\bsx \in \R^n : |\bsx-\bsu|=\sqrt{\lambda}\}. 
$$
It follows from the  equality $(2.2)$ that 
the inversion in $I(f)$ is written as   
$J_{I(f)}(\bsx)=\lambda J(\bsx-\bsu)+\bsu$. 
Thus the transformation $f(\bsx)=\lambda PJ(\bsx-\bsu)+\bsv$ 
can be  written as 
$$
f(\bsx)=P(J_{I(f)}(\bsx)-\bsu)+\bsv. 
$$
From this, one can observe that $f$ maps the interior of 
$I(f)=\{\bsx \in \R^n : |\bsx-\bsu|=\sqrt{\lambda}\}$ 
onto the exterior of  
$I(f^{-1})=\{\bsx \in \R^n : |\bsx-\bsv|=\sqrt{\lambda}\}$. 
Especially, the map $f$ restricted to $I(f)$ is an isometry. 

\subsection{Poincar\'{e} extensions.} 

We embed $\hR^n$ in $\hR^{n+1}$ 
in the natural way by  
$(x_1, \ldots, x_n) \in \R^n \mapsto (x_1, \ldots, x_n, 0) \in \R^{n+1}$ 
and $\fty \mapsto \fty$. 
Then for every $f \in \Mob(\hR^n)$, there is a unique 
$\tilde f \in \Mob(\hR^{n+1})$, called the {\it Poincar\'{e} extension} of $f$, 
such that $\tilde f|_{\hR^n}=f$. 
In fact,  if $f$ is a composition of inversions in 
$(n-1)$-dimensional spheres  
$\sigma_1, \ldots, \sigma_k$ in $\hR^n$, 
$\tilde f$ is obtained by the composition of inversions in 
 $n$-dimensional spheres 
$\tilde \sigma_1, \ldots, \tilde \sigma_k$ in $\hR^{n+1}$,   
where $\tilde \sigma_i$ is the sphere which 
is  orthogonal to $\hR^n$ at $\sigma_i$. 
In this way, $\Mob(\hR^n)$ can be regarded as a subgroup of 
$\Mob(\hR^{n+1})$. 
Let
$$
\HH^{n+1}=\{(x_1, \ldots ,x_{n+1}) \in \R^{n+1}:x_{n+1}>0\} 
$$
denote the upper half-space of $\hR^{n+1}$ 
equipped with the hyperbolic metric 
$ds^2=(dx_1^2+ \cdots +dx_{n+1}^2)/{x_{n+1}^2}$. 
It is the upper half-space model of the 
$(n+1)$-dimensional hyperbolic space. 
For every $f \in \Mob(\hR^n)$, 
its Poincar\'{e} extension  
$\tilde f \in \Mob(\hR^{n+1})$ induces an 
orientation preserving isomorphism of $\HH^{n+1}$. 
In this way, $\Mob(\hR^n)$ 
can be identified with 
the group $\Isom(\HH^{n+1})$ 
of orientation preserving isomorphism of $\HH^{n+1}$.  

\subsection{Classifications of M\"{o}bius transformations.}

\begin{defn}
 Let $f \in \Mob(\hR^n)=\Isom(\HH^{n+1})$. 
 We say that $f$  is {\it elliptic} if 
 $f$ has a fixed point in $\HH^{n+1}$. 
 If $f$ is not elliptic, $f$ is said to be  {\it parabolic} 
 if $f$ has exactly one fixed point in $\hR^n$, and  
 {\it loxodromic} if  
$f$ has exactly two fixed point in $\hR^n$. 
\end{defn}

It is known that every element $f \in \Mob(\hR^n)$ is 
either elliptic, parabolic  
or loxodromic (see (2.23) in \cite{Mat}). 
To describe a standard form of 
elliptic transformations, 
it is convenient to consider in 
the unit ball model 
${\bf  B}^{n+1}=\{\bsx \in \R^{n+1}:|\bsx|<1\}$ 
of the hyperbolic space: 
In this setting, if $f \in \Isom(\B^{n+1})$ is elliptic 
with $f(\bs0)=\bs0$, 
then $f(\bsx)=P(\bsx)$ for some $P \in SO(n+1)$. 
In the rest of the paper,  
we always consider in the upper half-space model of the hyperbolic space.  
Standard forms of parabolic and loxodromic transformations 
are given in the following: 
\begin{lem}[(2.24) in \cite{Mat}]
 Let $f \in \Mob(\hR^n)$. 
\begin{enumerate}
  \item   If $f$ is a parabolic transformation with a fixed point $\fty$, 
  then   $f(\bsx)=P(\bsx)+\bsu$ for some 
  $\bsu \in \R^n \sm \{\bs0\}$ and $P \in SO(n)$ with $P(\bsu)=\bsu$. 
  \item If $f$ is a loxodromic transformation with fixed points $\bs0, \fty$, 
  then $f(\bsx)=\lambda P(\bsx)$ for some 
$\lambda >0$ with $\lambda \ne 1$ and $P \in SO(n)$. 
\end{enumerate}
\end{lem}

\begin{defn}[pure parabolic, screw parabolic]
A parabolic transformation $f \in \Mob(\hR^n)$ is said to be 
{\it pure parabolic} if it is conjugate to a translation 
$\bsx \mapsto \bsx+\bsu$, and {\it screw parabolic} otherwise. 
\end{defn}

We now observe the case of $n=3$. 
Since every $P \in SO(3)$ has an eigenvalue of $1$, we have the following: 
\begin{lem}
Let $f \in \three$. 
\begin{enumerate}
  \item 
If $f$ is parabolic, it is conjugate in $\three$ to a transformation given by 
\begin{eqnarray}
\left(
  \begin{array}{c}
    x_1   \\
    x_2   \\
    x_3   \\
  \end{array}
\right)
 \mapsto \left(
  \begin{array}{ccc}
    \cos \theta   & -\sin \theta   & 0   \\
     \sin \theta  &  \cos \theta  & 0   \\
    0   &  0  &  1  \\
  \end{array}
\right)
\left(
  \begin{array}{c}
    x_1   \\
    x_2   \\
    x_3   \\
  \end{array}
\right)+
\left(
  \begin{array}{c}
    0   \\
    0   \\
    1   \\
  \end{array}
\right),  
\end{eqnarray}
where $0 \le \theta <2\pi$. 
( Then $f$ is pure parabolic if and only if $\theta=0$.)  
\item 
If $f$ is loxodromic, it is conjugate in $\three$ to a transformation given by 
\begin{eqnarray*}
\left(
  \begin{array}{c}
    x_1   \\
    x_2   \\
    x_3   \\
  \end{array}
\right)
 \mapsto \lambda \left(
  \begin{array}{ccc}
    \cos \theta   & -\sin \theta   & 0   \\
     \sin \theta  &  \cos \theta  & 0   \\
    0   &  0  &  1  \\
  \end{array}
\right)
\left(
  \begin{array}{c}
    x_1   \\
    x_2   \\
    x_3   \\
  \end{array}
\right), 
\end{eqnarray*}
where $\lambda>1$ and $0 \le \theta<2\pi$. 
\end{enumerate}
\end{lem}

It is worth to note that in $\three$ 
a pure parabolic transformation is topologically conjugate to 
any screw parabolic transformation. 
More precisely, 
a translation $(x_1,x_2,x_3) \mapsto (x_1,x_2,x_3)+(0,0,1)$ is conjugate to 
a screw parabolic transformation 
of the form $(2.3)$ 
by an orientation-preserving  homeomorphism $f:\hR^3 \to \hR^3$ defined by 
\begin{eqnarray*}
\left(
  \begin{array}{c}
    x_1   \\
    x_2   \\
    x_3   \\
  \end{array}
\right) 
\mapsto 
\left(
  \begin{array}{ccc}
    \cos (x_3 \theta)   & -\sin (x_3 \theta)   & 0   \\
     \sin (x_3 \theta)   &  \cos (x_3 \theta)  & 0   \\
    0   &  0  &  1  \\
  \end{array}
\right)\left(
  \begin{array}{c}
    x_1   \\
    x_2   \\
    x_3   \\
  \end{array}
\right). 
\end{eqnarray*}

\subsection{Kleinian groups.}
A discrete subgroup 
$\Gamma$ of $\Isom(\HH^{n+1})=\Mob(\R^n)$ 
is called a $(n+1)$-dimensional {\it Kleinian group}. 
It is known that $\Gamma$ is a Kleinian group 
if and only if $\Gamma$ acts  properly discontinuously on $\HH^{n+1}$; 
that is,  given any compact subset $K \subset \HH^{n+1}$, 
the set $\{\gamma  \in \Gamma : \gamma(K) \cap K \ne \emptyset\}$ is finite.  

Let $\Gamma$ be a $(n+1)$-dimensional Kleinian group.   
The {\it domain of discontinuity} $\Omega(\Gamma) \subset \hR^n$ of $\Gamma$ 
is defined to be  the set of points 
$\bsx \in \hR^n$ such that there is a neighborhood $U$ of 
$\bsx$ such that the set 
$\{\gamma \in \Gamma : \gamma(U) \cap U= \emptyset\}$ is finite. 
The complement of $\Omega(\Gamma)$ in $\hR^n$ 
is called the {\it limit set} of $\Gamma$, 
and denoted by  $\Lambda(\Gamma)$.  

\subsection{Notations.}
Throughout this paper, we use the following notations: 
\begin{itemize}
  \item If $S$ is a subset of a group $G$,  
we denote by $\la S \ra \subset G$ the subgroup generated by $S$. 
\item According to the usual convention, 
a M\"{o}bius transformation $f \in \two$ is written 
in a linear fractional form 
$$f(\tau)=\frac{a\tau+b}{c\tau+d} \quad (a,b,c,d \in \C) 
$$
by identifying $\hR^2$ with the Rieman sphere $\hC$. 
\item A M\"{o}bius transformation $f:\hR^3 \to \hR^3$ in $\three$ is 
written as a function $f(\bsx)$ of $\bsx=(x,y,z) \in \R^3$.  
\item
We denote by $P_{z=0}$ the plane 
$\{(x,y,z) \in \R^3:z=0\}$ and 
by $P_{r=0}$ the plane 
$\{(p,q,r) \in \R^3:r=0\}$ and so on.  
In addition, we write $\hP_{z=0}=P_{z=0} \cup \{\fty\}$  
and $\hP_{r=0}=P_{r=0} \cup \{\fty\}$ and so on. 
\end{itemize}

\section{Groups in $\two$} 

In this section, we recall some basic facts of 
Kleinian groups in $\two$ which are isomorphic 
to the fundamental groups of a thrice-punctured sphere,   
a fourth-punctured sphere or a punctured torus. 
We refer the reader to Kra \cite{Kr} and Keen--Series \cite{KS}. 

A surface $\Si_{g,n}$ of {\it type} $(g,n)$  
 is an oriented closed surface of genus $g$ 
with $n$ punctures. 
If $n > 0$, the fundamental group 
$\pi_1(\Si_{g,n})$ of $\Si_{g,n}$ is a free group. 
A representation $\rho:\pi_1(\Si_{g,n}) \to \two$ is said to be 
{\it type-preserving} if it takes 
a loop surrounding a puncture to a parabolic transformation.  
In this section, 
we consider images of faithful type-preserving representations 
$\rho:\pi_1(\Si_{g,n}) \to \two$ 
for 
$(g,n)=(0,3), \, (0,4)$ and $(1,1)$. 
We choose generators for 
$\pi_1(\Si_{0,3})=\la b,c \ra$, 
$\pi_1(\Si_{0,4})=\la b,c,d \ra$ and 
$\pi_1(\Si_{1,1})=\la a,b \ra$ as in Figure 3. 
\begin{figure}[htbp]
  \begin{center}
\includegraphics[keepaspectratio=true,height=30mm]{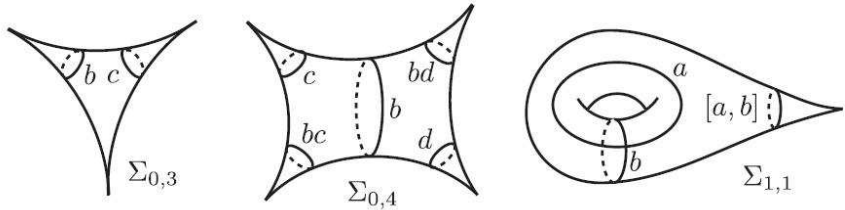}
  \end{center}
  \caption{Surfaces $\Si_{0,3}$, $\Si_{0,4}$ and $\Si_{1,1}$,  and 
free homotopy classes 
of elements of their fundamental groups. }
\end{figure}
We define groups of type $(0,3),\,(0,4)$ and $(1,1)$ in $\two$ as follows: 
\begin{defn}
Let $\Gamma$  be a subgroup of $\two$. 
\begin{itemize}
  \item $\Gamma$ is said to be  of {\it type} $(0,3)$  
  if  it is the image of a 
  faithful, type-preserving representation 
   $\rho:\pi_1(\Si_{0,3})=\la b,c \ra \to \two$, 
   or equivalently, if it is a rank-$2$ free group $\la \beta,\gamma \ra$ such that 
   $\beta, \gamma$ and $\beta\gamma$ are parabolic. 
   \item $\Gamma$ is said to be of {\it type} $(0,4)$  
   if  it is the image of a
  faithful, type-preserving representation 
   $\rho:\pi_1(\Si_{0,4})=\la b,c,d \ra \to \two$ 
   such that $\rho(b)$ is parabolic, or equivalently, 
if it is a rank-$3$ free group $\la \beta,\gamma, \delta \ra$ 
such that $\beta,\gamma,\delta, \beta\gamma$ and $\beta\delta$ are parabolic.  
   \item $\Gamma$ is said to be of {\it type} $(1,1)$  
   if it is the image of  a  
  faithful, type-preserving representation 
   $\rho:\pi_1(\Si_{1,1})=\la a,b \ra \to \two$ 
   such that $\rho(b)$ is parabolic, or equivalently, 
   if it is a rank-$2$ free group $\la \alpha, \beta \ra$ 
      such that $\beta$ and $[\alpha,\beta]$ is parabolic. 
\end{itemize}
\end{defn}

\subsection{Groups of type $(0,3)$}

Observe  that the subgroup $\la B,C \ra$ of $\two$ generated by 
\begin{eqnarray*}
B(\tau)=\tau+2 \quad \text{and} \quad 
C(\tau)=\frac{1}{B(\frac{1}{\tau})}=\frac{\tau}{2\tau+1}
\end{eqnarray*}
is of type $(0,3)$. 
It is known that any group of type $(0,3)$ in $\two$ is conjugate to this group: 
\begin{lem}
Let 
$\Gamma=\la \beta,\gamma \ra$ be a rank-$2$ free subgroup of $\two$ 
such that $\beta, \gamma$ and $\beta \gamma$ are parabolic. 
Then $\Gamma$ is conjugate in $\two$ 
to the group $\la B,C\ra$ defined above. 
\end{lem}

\subsection{Groups of type $(0,4)$}

For a given $\mu \in \C$, let  
$$
H_\mu=\la B,C,D_\mu\ra
$$ 
be the group in $\two$ 
generated by 
\begin{eqnarray*}
B(\tau)=\tau+2,  \quad  
C(\tau)=\frac{\tau}{2\tau+1} \quad \text{and} \quad  
D_\mu(\tau)=C(\tau-\mu)+\mu. 
\end{eqnarray*}
It is known that any group of type 
$(0,4)$ in $\two$ is normalized in this form: 

\begin{lem}
Let $\Gamma=\la \beta,\gamma,\delta \ra$ be a rank-$3$ free subgroup of $\two$ 
such that $\beta,\gamma,\delta, \beta\gamma$ and $\beta \delta$ are parabolic. 
Then $\Gamma$ is conjugate in $\two$ 
to $H_\mu=\la B,C,D_\mu \ra$ for some $\mu \in \C$. 
\end{lem}
A fundamental domain for $H_\mu$ with $\im \mu \ge 1$ 
is given in the left of Figure 4.  
The {\it Maskit slice of groups of type $(0,4)$ in $\two$} is defined by 
\begin{eqnarray*}
\M_{0,4}=\{\mu \in \C: 
H_\mu=\la B,C,D_\mu \ra \text{ is a rank-$3$ free Kleinian group}\}. 
\end{eqnarray*}
Some basic properties of $\MN$ can be found in the next subsection. 
\begin{figure}[htbp]
  \begin{center}
    \includegraphics[keepaspectratio=true,height=70mm]{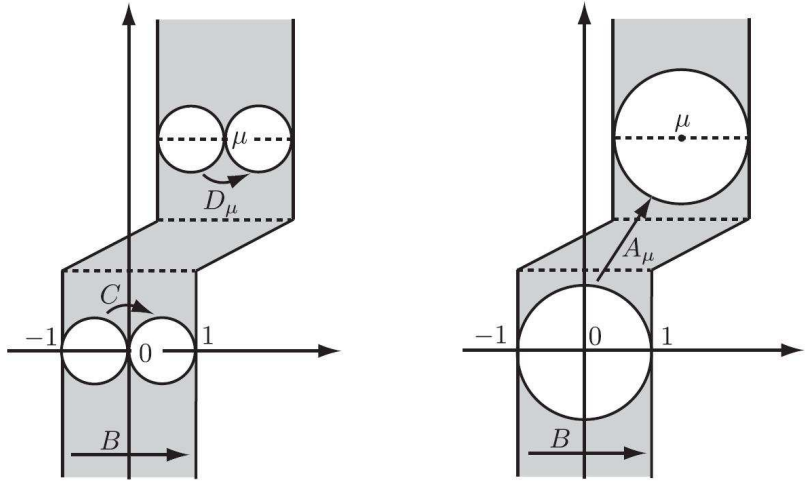}
    \end{center}
  \caption{Fundamental domains (the shaded regions) 
  for $H_\mu$ (left) and $G_\mu$ (right). }
\end{figure}

\subsection{Groups of type $(1,1)$}

For a given $\mu \in \C$, let  
$$
G_\mu=\la A_\mu ,B \ra
$$ 
be the group in $\two$ generated by 
\begin{eqnarray*}
A_\mu(\tau)=\frac{1}{\tau}+\mu \quad \text{and} \quad 
B(\tau)=\tau+2. 
\end{eqnarray*}
It is known that any group of type $(1,1)$ in $\two$ is normalized in this form: 

\begin{lem}
Let $\Gamma=\la \alpha,\beta \ra$ be a rank-$2$ free subgroup of $\two$ 
such that $\beta$ and $[\alpha, \beta]$ are parabolic. 
Then $\Gamma$ is conjugate in $\two$ 
to $G_\mu=\la A_\mu, B \ra$  for some $\mu \in \C$. 
\end{lem}

A fundamental domain for $G_\mu$ with $\im \mu \ge 2$ is given 
in the right of Figure 4.  
The {\it Maskit slice of groups of type 
$(1,1)$ in $\two$} is defined by 
\begin{eqnarray*}
\MM=\{\mu \in \C: 
G_\mu=\la A_\mu, B \ra  \text{ is a rank-$2$ free Kleinian group}\} 
\end{eqnarray*}
(see also Figure 5).  
It is easy to see that $\MM$ is invariant under the maps 
$\mu \mapsto \mu+2$ and $\mu \mapsto -\mu$, and that 
$\MM$ contains the set $\{\mu \in \C : |\im \mu| \ge 2 \}$. 
It is also known that the Maskit slice $\MM$ is contained in 
the set $\{\mu \in \C: |\im \mu|> 1\}$ (cf. \cite{KS}). 
Furthermore, Minsky \cite{Mi} showed that 
$\MM $consists of two connected component and that 
each connected component of 
$\MM$ plus $\{\fty\}$ is a closed topological disc in $\hat \C$.

We mention here the relationship between the Maskit slices $\MN$ and $\MM$. 
First observe for every $\mu \in \C$ that 
$$
C=A_\mu^{-1}BA_\mu \quad  \text{and} \quad D_\mu=A_\mu BA_{\mu}^{-1}. 
$$
Thus $H_\mu=\la B,C,D_\mu\ra$ is always a subgroup of 
$G_\mu=\la A_\mu,B \ra$. 
This implies that $\MM \subset \M_{0,4}$. 
Furthermore, it is known by Kra that $\MN$ is similar to $\MM$:  
\begin{thm}[Kra \cite{Kr}]
The map $\tau \mapsto 2\tau$ of $\C$ onto itself induces 
a bijective map from $\MN$ onto $\MM$.  
\end{thm}
\begin{figure}[htbp]
  \begin{center}
  \includegraphics[keepaspectratio=true,height=50mm]{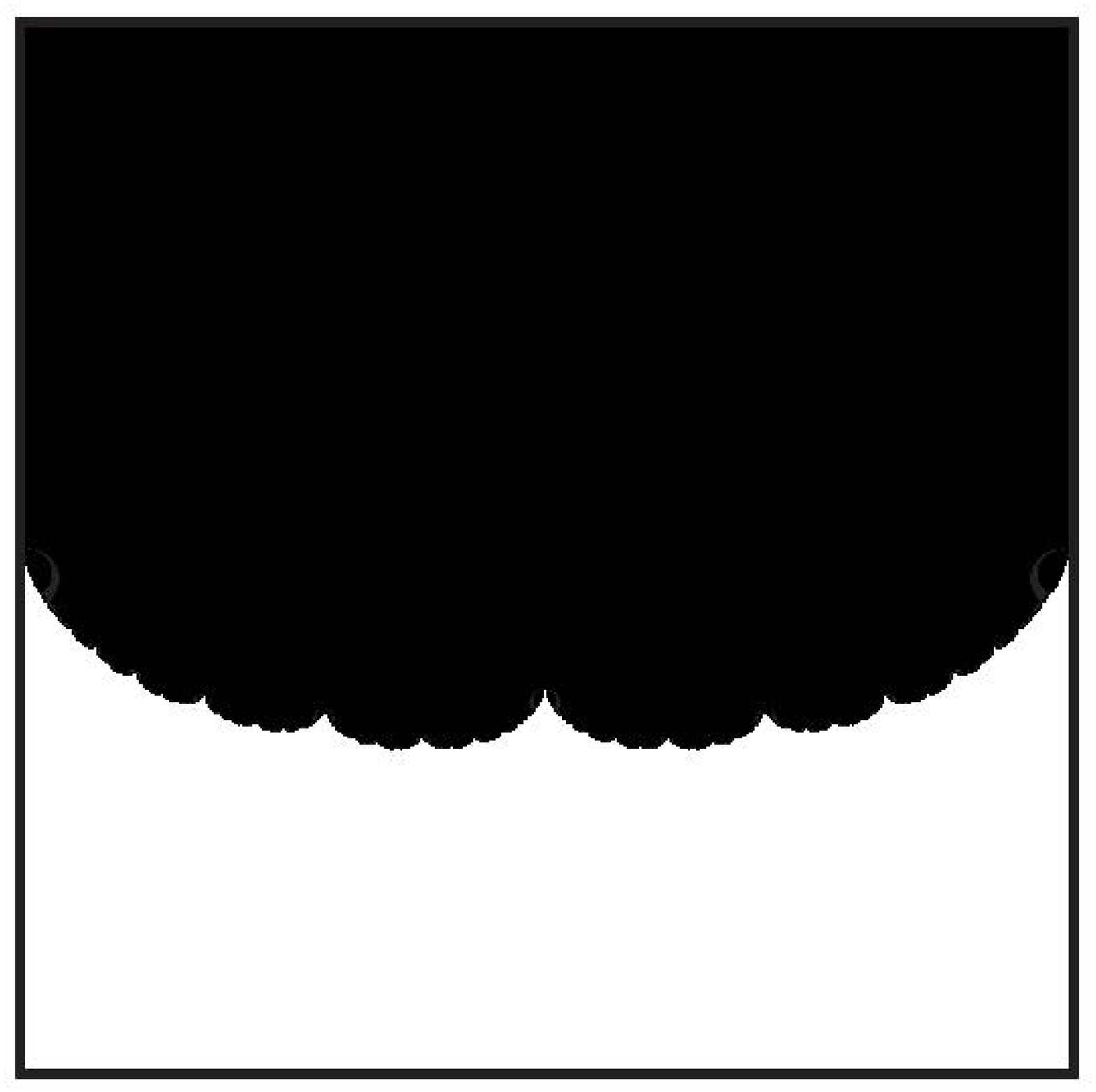}
  \end{center}
  \caption{The Maskit slice $\MM$ (black part) restricted to the domain  
  $\{\mu \in \C : 0 \le \re \mu \le 2, \, 1 \le \im \mu \le 3\}$.}
\end{figure}

\section{Groups in $\three$} 

We define $4$-dimensional analogues of $3$-dimensional 
groups of types $(0,3),\,(0,4)$ and $(1,1)$ as follows: 
\begin{defn}
Let $\Gamma$  be a subgroup of $\three$. 
\begin{itemize}
  \item $\Gamma$ is said to be  of {\it type} $(0,3)$  
  if it is a rank-$2$ free group $\la \beta,\gamma \ra$ such that 
   $\beta, \gamma$ and $\beta\gamma$ are pure parabolic. 
   \item $\Gamma$ is said to be of {\it type} $(0,4)$  
if it is a rank-$3$ free group $\la \beta,\gamma, \delta \ra$ 
such that $\beta,\gamma,\delta, \beta\gamma$ and $\beta\delta$ are pure parabolic.  
   \item $\Gamma$ is said to be of {\it type} $(1,1)$  
   if it is a rank-$2$ free group $\la \alpha, \beta \ra$ 
      such that $\beta$ and $[\alpha,\beta]$ are pure parabolic. 
\end{itemize}
\end{defn}
As mentioned in Introduction, 
it is just a technical reason why we assume some elements of these groups 
are (not only parabolic but also) pure parabolic. 
In this section, we study normalizations of these groups in $\three$ up to conjugation. 

\subsection{Groups of type $(0,3)$ in $\three$}

Observe that the Poincar\'{e} extension of the group 
$\la B(\tau)=\tau+2, C(\tau)=\tau/(2\tau+1) \ra \subset \two$ 
acting on $\hC$ is given by 
\begin{eqnarray*}
\la B(\bsx)=\bsx+(2,0,0), \ 
C(\bsx)=JBJ(\bsx) \ra \subset \three, 
\end{eqnarray*}
where $\hC$ is identified with $\hP_{z=0} \subset \hR^3$ via the map 
$\tau=x+iy \mapsto \bsx=(x,y,0)$, and 
$J(\bsx)=\bsx/|\bsx|^2$ is the inversion in the unit sphere. 
(Here and hereafter,  we use the same symbols $B, C$ for 
the Poincar\'{e} extensions of $B,C \in \two$ by abuse of notation.) 
Any group of type $(0,3)$  in $\three$ is conjugate to this group: 

\begin{thm}
Let $\Gamma=\la \beta,\gamma \ra$ be a rank-$2$ free subgroup of $\three$ such that 
$\beta$, $\gamma$ and $\beta \gamma$ are pure parabolic. 
Then $\Gamma$ is conjugate in $\three$ 
to the group $\la B, C \ra$ defined above.  
\end{thm}

\begin{proof}
Observe first that $\Fix(\beta) \ne \Fix(\gamma)$.  
In fact, if $\Fix(\beta)=\Fix(\gamma)$, 
$\beta$ commutes with $\gamma$ 
since both $\beta,\gamma$ are pure parabolic. 
This contradicts that  
$\la \beta, \gamma \ra$ is free. 
Therefore, after conjugating the group 
$\la \beta,\gamma \ra$ in $\three$ if necessary, 
we may assume that 
$\beta(\fty)=\fty$ and $\gamma(\bs0)=\bs0$. 
Observe that there exist lines $l$ and $m$ through the origin $\bs0$ 
which are invariant under the actions of $\la \beta \ra$ and 
$\la \gamma \ra$, respectively. 
After conjugating by an orthonormal translation if necessary, 
we may assume  that both $l,\,m$ lie 
in the sphere $\hP_{z=0}$. 
Then one can see that the group 
$\la \beta, \gamma \ra$ preserves the sphere 
$\hP_{z=0}$ and its orientation. 
Thus the result follows form Lemma 3.2. 
\end{proof}

\begin{rem}
Recently, 
Y.~Kim and K.~Sakugawa independently announced (oral communication) that there is a continuous family 
of rank-$2$ free Kleinian groups $\Gamma=\la \beta, \gamma \ra$ in $\three$ such that 
$\beta, \gamma$ and $\beta\gamma$ are (screw) parabolic.  
\end{rem}

In the next subsection, we need the following: 

\begin{cor}
Let $\la \beta,\gamma \ra$ be a rank-$2$ free subgroup of $\three$ such that 
$\beta$, $\gamma$ and $\beta \gamma$ are pure parabolic. 
If $\beta=B$ and $\gamma(\bs0)=\bs0$ then $\gamma=C^{\pm 1}$, that is, 
$\la \beta, \gamma \ra=\la B,C \ra$. 
\end{cor}

\begin{proof}
By Theorem 4.2, there is a transformation $f \in \three$ 
which conjugates $\la B, \gamma \ra$ to $\la B,C \ra$.   
Since $J$ conjugates $B$ and $C$ each other, 
we may assume that $fBf^{-1}=B$ and $f \gamma f^{-1}=C^{\pm 1}$. 
This implies that $f(\bs0)=\bs0$ and $f(\fty)=\fty$. 
Then by Lemma 2.1, we have 
$f=\lambda P$, where $\lambda>0$ and $P \in SO(3)$. 
Moreover, it follows from the condition 
$fBf^{-1}=B$ that $f$ is a rotation along the $x$-axis. 
Since such $f$ conjugates $C$ to itself, 
we have $\gamma=f^{-1}C^{\pm 1} f=C^{\pm 1}$.  
\end{proof}

\subsection{Groups of type $(0,4)$}

In this subsection, we will show that 
every group of type $(0,4)$ in $\three$ is 
conjugate to a group of type $(0,4)$ in $\two$. 

For a given $\bsp=(p,q,r) \in \R^3$, let
\begin{eqnarray*}
H_\bsp=\la B,C,D_\bsp \ra
\end{eqnarray*}
be the subgroup of $\three$ generated by 
\begin{eqnarray*}
B(\bsx)=\bsx+(2,0,0), \quad
C(\bsx)=JBJ(\bsx) \quad  \text{and} \quad  D_\bsp(\bsx)=C(\bsx-\bsp)+\bsp. 
\end{eqnarray*}
Observe that if $\bsp$ lies in the plane $P_{r=0}$,  
the group $H_\bsp$ for $\bsp=(p,q,0)$ is the Poincar\'{e} extension of the group 
$H_\mu=\la B,C,D_\mu \ra$ for $\mu=p+iq$ defined in 3.2. 

\begin{thm}
Let $\Gamma=\la \beta,\gamma,\delta \ra$ be a rank-$3$ free subgroup of $\three$ such that 
$\beta, \gamma, \delta, \beta\gamma$ and $\beta\delta$ are pure parabolic. 
Then $\Gamma$ is conjugate in $\three$ to $H_\bsp=\la B, C, D_\bsp \ra$ 
for some $\bsp=(p,q,0) \in P_{r=0}$. 
\end{thm}

\begin{proof}
Note that the subgroups $\la \beta, \gamma \ra$, $\la \beta, \delta  \ra$ of 
$\la \beta, \gamma, \delta \ra$ are of type $(0,3)$ in $\three$. 
Therefore, after conjugating 
$\la \beta, \gamma, \delta \ra$ in $\three$ if necessary, 
we may assume that $\beta=B$ and $\Fix(\gamma)=\bs0$, and hence that 
$\la \beta, \gamma \ra=\la B,C\ra$  by Corollary 4.3. 
Since a rotation along the $x$-axis 
conjugates the group $\la B,C \ra$ to itself, 
we may also assume that the fixed point $\bsp:=\Fix(\delta)$ of $\delta$ 
lies in the plane $P_{z=0}$. 
To obtain the result, it suffices to show that 
$\la \beta,\delta \ra=\la B,D_\bsp \ra$. 
Now let us conjugate the group $\la B, \delta \ra$ by 
the translation $T_\bsp^{-1}(\bsx)=\bsx-\bsp$
to the group $T_\bsp^{-1} \la B,  \delta \ra T_\bsp =\la B, T_\bsp^{-1} \delta T_\bsp \ra$. 
 Since $T_\bsp^{-1} \delta T_\bsp$ fixes $\bs0$, we have 
$\la B, T_\bsp^{-1} \delta T_\bsp \ra=\la B,C \ra$ from Corollary 4.3.  
Thus we conclude that $\la B,\delta \ra=T_\bsp \la B,C \ra T_\bsp^{-1}=\la B,D_\bsp \ra$. 
\end{proof}

\subsection{Groups of type $(1,1)$}

In this subsection,  
we will obtain a normalization of groups of type 
$(1,1)$ in $\three$.  
In contrast to the case of groups of type $(0,4)$, 
we will see that 
the space of groups of type $(1,1)$ in $\three$ 
is strictly larger than the space of groups of type $(1,1)$ in $\two$. 

For a given $\bsp \in \R^3$, let 
\begin{eqnarray*}
G_\bsp=\la A_\bsp, B \ra
\end{eqnarray*}
be the subgroup of $\three$ generated by 
\begin{eqnarray*}
A_\bsp(\bsx)=\hat J J(\bsx)+\bsp  \quad \text{and} \quad 
B(\bsx)=\bsx+(2,0,0), 
\end{eqnarray*}
where  
$$
\hat J(x,y,z)=(x,-y,z)  
$$
is the inversion in the sphere $\hP_{y=0}$. 
Observe for every $\bsp \in \R^3$ that 
\begin{eqnarray*}
C=A_\bsp^{-1}BA_\bsp \quad \text{and} \quad D_\bsp=A_\bsp B A_\bsp^{-1}. 
\end{eqnarray*}
Thus $H_\bsp=\la B,C,D_\bsp \ra$ is always a subgroup of 
$G_\bsp=\la A_\bsp, B \ra$. 
In addition, observe that if $\bsp$ lies in the plane $P_{r=0}$,  
the group $G_\bsp$ for $\bsp=(p,q,0)$ is the Poincar\'{e} extension of the group 
$G_\mu=\la A_\mu, B\ra$ for $\mu=p+iq$ defined in 3.3. 
The next theorem reveals that 
every group of type $(1,1)$ in $\three$
is conjugate to the group $G_\bsp$ for some $\bsp \in \R^3$. 
(We will prove a slightly stronger statement.)  

\begin{thm}
Let 
$\Gamma=\la \alpha, \beta \ra$ be a rank-$2$ free group in $\three$ such that 
$\beta$ is pure parabolic and $[\alpha,\beta]$ is (not necessarily pure) parabolic. 
Then $\Gamma$ is conjugate in $\three$ to $G_\bsp=\la A_\bsp,B \ra$ 
for some $\bsp=(p,q,r) \in \R^3$. 
\end{thm}

\begin{proof}
We first show that $\Fix(\alpha) \cap \Fix(\beta)=\emptyset$. 
If not, we may assume that $\fty \in \Fix(\alpha) \cap \Fix(\beta)$. 
It then follows from Lemma 2.1 that  
$\alpha(\bsx)=\lambda P(\bsx)+\bsu$ and $\beta(\bsx)=\bsx+\bsv$ 
for some $\lambda>0$, $P \in SO(3)$ and $\bsu,\,\bsv \in \R^3$. 
Then a calculation yields $[\alpha,\beta](\bsx)=\bsx+\lambda P(\bsv)-\bsv$. 
Thus $[\alpha,\beta]$ is pure parabolic with the fixed point $\fty$. 
This implies that $\beta$ commutes with $[\alpha,\beta]$, which contradicts that 
$\la \alpha,\beta \ra$ is free. 

Therefore we may assume that  
$\beta(\infty)=\infty$ and $\alpha(\bs 0)=\infty$.  
It then follows from Lemma 2.1 that 
$\alpha$, $\beta$ are of the forms 
$\alpha({\bs x})=\lambda PJ({\bs x})+{\bs p}$ and 
$\beta({\bs x})={\bs x}+{\bs u}$ for some 
$\lambda>0$, $P \in O(3) \sm SO(3)$ and 
${\bs p}, \,{\bs u} \in \R^3$. 
We may further conjugate the group $\la \alpha, \beta \ra$ by a 
transformation $f \in \three$ fixing $\bs0$ and $\fty$ 
without changing our assumptions $\beta(\infty)=\infty$, $\alpha(\bs 0)=\infty$.  
Such a transformation is of the form 
$f=\lambda' Q \in \three$ with $\lambda'>0,\,Q \in SO(3)$ by Lemma 2.3.   
Therefore, after conjugating $\la \alpha,\beta \ra$ 
by a magnification if necessary,  
we may assume that the radius of the isometric sphere of 
$\alpha$ equals $1$, and hence that 
$\alpha,\,\beta$ are of the forms  
$\alpha({\bs x})=PJ({\bs x})+{\bs p}$ and 
$\beta({\bs x})={\bs x}+{\bs u}$. 
In addition, we claim that we may also assume that 
$$
{\bs u}=(u,v,0) \quad  \text{and} \quad P^{-1}(\bs u)=(u,-v,0)
$$
for some $u,v \ge 0$ after conjugating $\la \alpha,\beta \ra$ by a 
suitable orthonormal transformation if necessary. 
In fact, choose $Q \in SO(3)$ so that 
$Q(\bsu)=(u,v,0)$ and $QP^{-1}(\bsu)=(u,-v,0)$. 
Then $Q$ conjugates $\alpha(\bsx)=PJ(\bsx)+\bsp$ to 
$Q \alpha Q^{-1}(\bsx)=QPQ^{-1}J(\bsx)+Q(\bsp)$ and $\beta(\bsx)=\bsx+\bsu$ to 
$Q \beta Q^{-1}(\bsx)=\bsx+Q(\bsu)$.  
Now observe that $(QPQ^{-1})^{-1}$ takes $Q(\bsu)=(u,v,0)$ to 
$QP^{-1}(\bsu)=(u,-v,0)$.  
Therefore, replacing $QPQ^{-1}$, $Q(\bsp)$ and $Q(\bsu)$ by 
$P$, $\bsp$ and $\bsu$, respectively, 
we obtain the claim. 

Next we will show that $\bsu=P^{-1}(\bsu)=(2,0,0)$ 
by using the condition 
that $[\alpha,\beta]$ is parabolic. 
Recall that $T_{\bs v}$ denotes the translation 
$\bsx \mapsto {\bs x} +{\bs v}$ with ${\bs v} \in \R^3$. 
Using this notation, we have   
$\alpha({\bs x})=T_{\bs p} P J ({\bs x}), \, 
\beta({\bs x})=T_{\bs u}(\bs x), \,  
\alpha^{-1}=J P^{-1} T_{-{\bs p}}(\bs x)$ and 
$\beta^{-1}(\bs x)=T_{-{\bs u}}(\bs x)$. 
Then  the transformation  
$(\beta\alpha)^{-1}[\alpha,\beta]\beta\alpha=\alpha^{-1}\beta^{-1}\alpha\beta$ 
is calculated as 
\begin{eqnarray*}
\alpha^{-1}\beta^{-1}\alpha\beta(\bs x)
&=&J P^{-1}T_{-{\bs p}}T_{-{\bs u}}T_{\bs p} PJ T_{\bs u}(\bs x) \\
&=&J P^{-1}T_{-{\bs u}} P J T_{\bs u}(\bs x) \\
&=&J T_{-P^{-1}({\bs u})} J T_{\bs u}(\bs x), 
\end{eqnarray*}
where the last equality follows from 
$P^{-1}T_{-{\bs u}} P(\bsx)=P^{-1}(P(\bsx)-\bsu)=\bsx-P^{-1}(\bsu)=
T_{-P^{-1}(\bsu)}(\bsx)$. 
Since $\bsu=(u,v,0)$ and $P^{-1}(\bsu)=(u,-v,0)$ lie in the plane $P_{z=0}$, 
the transformation 
$\alpha^{-1}\beta^{-1}\alpha\beta=J T_{-P^{-1}({\bs u})} J T_{\bs u}$ 
preserves the sphere 
$\hP_{z=0}$ and its orientation. 
By letting $\mu=u+iv$ and identifying 
$\hP_{z=0}$ with $\hC$ via the map 
${\bs x}=(x,y,0) \mapsto \tau=x+iy$, 
we see that the action of 
$\alpha^{-1}\beta^{-1}\alpha\beta$  
restricted to the sphere $\hP_{z=0} \cong \hC$ 
is a M\"{o}bius transformation 
\begin{eqnarray*}
\tau \mapsto \dfrac{\tau+\mu}{-\mu \tau +1-\mu^2},  
\end{eqnarray*}
whose matrix representation is 
$$
\left(
  \begin{array}{cc}
    1  & \mu   \\
     -\mu  &  1-\mu^2  \\
  \end{array}
\right) \in \psl. 
$$
Since 
$\alpha^{-1}\beta^{-1}\alpha\beta$ is 
parabolic and $\re \mu =u \ge 0$, 
we conclude that  $\mu=2$,  
and hence that ${\bs u}=P^{-1}(\bsu)=(2,0,0)$. 

Now let us denote by $R_\theta \in SO(3)$ the rotation of angle $\theta \in \R$ 
along the $x$-axis. 
Since $P \in O(3) \sm SO(3)$ and $P(\bsu)=\bsu$, 
we have $P=R_\varphi \hat J$ for some $\varphi \in \R$.  
Noting that  
 $\hat J R_{-\theta}=R_{\theta} \hat J$ 
and $J R_{\theta}=R_{\theta} J$ for every $\theta \in \R$, 
we obtain 
\begin{eqnarray*}
&&R_\theta \alpha R_\theta^{-1}(\bsx)
=R_\theta (R_\varphi \hat J J R_{-\theta}(\bsx)+\bsp)
=R_{\varphi+2\theta}\hat J J(\bsx)+R_\theta (\bsp),  \\
&&R_\theta \beta R_\theta^{-1}(\bsx)=\beta(\bsx). 
\end{eqnarray*}
Therefore, after conjugating $\la \alpha,\beta \ra$ by $R_\theta$ with 
$\theta \equiv -\varphi/2$ $(\mathrm{mod} \ \pi)$ if necessary, 
we may assume that $\alpha, \beta$ are of the forms  
 $\alpha(\bsx)=A_\bsp(\bs x)=\hat JJ(\bs x)+{\bs p}$ and 
$\beta(\bsx)=B(\bs x)={\bs x}+(2,0,0)$. 
Thus we obtain the desired normalization. 
\end{proof}

\begin{rem}
The condition that the group $\Gamma=\la \alpha,\beta \ra$ 
is free is only used to show that $\Fix(\alpha) \cap \Fix(\beta)=\emptyset$. 
\end{rem}

\section{The space of groups of type $(1,1)$ in $\three$}

We define the {\it space of Kleinian groups of type $(1,1)$ in $\three$} by 
$$
\hMM=\{\bsp=(p,q,r) \in \R^3 : 
G_\bsp=\la A_\bsp, B \ra \text{ is a rank-$2$ free Kleinian group}\}. 
$$
The aim of this section is to study the shape of $\hMM \subset \R^3$. 
Especially, the slices of $\hMM$ along the planes 
$P_{q=0}$ and $P_{p=0}$ will be studied. 

\begin{rem}
Although $B, [A_\bsp,B] \in G_\bsp$ are pure parabolic for every $\bsp$, 
it remains the  possibility that for some  $\bsp \in \hMM$, the group 
$G_\bsp$  has (accidental) screw parabolic transformations. 
We do not know whether it does happen or not. 
\end{rem}

Observe that  by definition we have 
$$
\hMM \cap P_{r=0}=\MM, 
$$
where the set 
$\MM \subset \C$ is regarded as a subset of the plane $P_{r=0}$ 
via the identification 
$\C \ni p+iq \mapsto (p,q,0) \in P_{r=0}$. 

We next observe in Theorem 5.1 below that 
there are two subsets $V_1, V_2$ of the parameter space $\R^3=\{\bsp=(p,q,r)\}$ 
such that $V_1 \subset \hMM \subset V_2$ and that 
both of them are invariant under the action of the rotation along the $p$-axis. 
Recall that 
$R_\theta \in SO(3) \subset \three$ denotes the rotation of angle 
$\theta \in \R$ around the $x$-axis. 
We also regard $R_\theta$ acts on the parameter space 
$\R^3=\{\bsp=(p,q,r)\}$ as a rotation around the $p$-axis. 
Furthermore, identifying $\hC$ with $\hP_{r=0} \subset \hR^3$ via the map 
$\mu=p+iq \mapsto \bsp=(p,q,0)$, we define a map $R_\theta:\hC \to \hR^3$ by 
\begin{eqnarray*}
R_\theta(\mu):=R_\theta(\bsp)=(p,q \cos \theta, q \sin \theta). 
\end{eqnarray*}

\begin{thm}
We have  
\begin{eqnarray*}
\bigsqcup_{0 \le \theta <\pi}R_\theta(\{\mu \in \C :| \im \mu| \ge 2 \}) 
\subset \hMM \subset 
\bigsqcup_{0 \le \theta<\pi}R_\theta(\MN). 
\end{eqnarray*}
\end{thm}

\begin{proof}
Suppose first that $\bsp \in \hMM$. 
Then $G_\bsp$ is discrete, and thus  
the subgroup $H_\bsp \subset G_\bsp$ is also discrete. 
Choose $0 \le \theta <\pi$ and $\mu \in  \C$ such that 
$\bsp=R_\theta(\mu)$. 
Then we have 
$H_\bsp=R_\theta H_\mu R_\theta^{-1}$,  
where $H_\mu$ is naturally regarded as a subgroup of $\three$. 
It then follows from the discreteness of $H_\bsp$ that $\mu \in \MN$, and hence that  
$\bsp=R_\theta(\mu) \in \bigsqcup_{0 \le \theta<\pi}R_\theta(\MN)$. 

Next suppose that 
$\bsp \in R_\theta (\{\mu \in \C :| \im \mu| \ge 2 \})$ 
for some $0 \le \theta <\pi$. 
Then the interiors of the isometric spheres $I(A_\bsp)=\{\bsx \in \R^3 : |\bsx|=1 \}$ of 
$A_\bsp$ and 
$I(A_\bsp^{-1})=\{\bsx \in \R^3 : |\bsx-\bsp|=1 \}$ of 
$A_\bsp^{-1}$ are disjoint. 
Therefore we can find a fundamental domain for $G_\bsp=\la A_\bsp, B \ra$ 
in $\hR^3$ as a $3$-dimensional analogue of the fundamental domain for $G_\mu$ 
in $\hC$ as in Figure 4. 
Then by Poincar\'{e}'s polyhedron theorem 
(see Maskit \cite[IV, H]{Mas} and the remark below), 
$G_\bsp$ is a rank-$2$ free Kleinian group in $\three$. 
Thus we conclude that $\bsp \in \hMM$. 
\end{proof}

\begin{rem}
Poincar\'{e}'s polyhedron theorem for 
$3$-dimensional polyhedra in $\hR^3$ with 
side-pairing maps in $\three$ can be deduced from 
the theorem for 
$4$-dimensional polyhedra in $\HH^4$ with 
side-pairing maps in $\Isom (\HH^4)$; the precise statement 
can be found in Maskit \cite{Mas}. 
We remark that, in general, a polyhedron in $\hR^3$  need not be convex. 
See also Epstein--Petronio \cite{EP}. 
\end{rem}

\subsection{Slice through the plane $P_{q=0}$}

The goal of this subsection is to show Theorem 5.4, which states that 
the Maskit slice $\MN$ of groups of type $(0,4)$ in $\two$ 
appears as the slice of $\hMM$ through the plane 
$P_{q=0}=R_{\frac{\pi}{2}}(\C)$.  

Before starting the proof of Theorem 5.4,  
we will make some observations on groups  
$G_\bsp \subset \three$ with $\bsp \in P_{q=0}$.  
Observe that if $\bsp \in P_{q=0}$, the action of 
$G_\bsp=\la A_\bsp, B\ra$ preserves 
the sphere $\hP_{y=0}$. 
More precisely, 
let $\bsp=(p,0,r) \in P_{q=0}$ and $\mu=p+ir \in \C$. 
Then the action of $G_\bsp=\la A_\bsp, B \ra$ 
restricted to the sphere $\hP_{y=0}$ is given by 
$$
\ch G_\mu=\langle \ch A_\mu, B \rangle; \quad 
\ch A_\mu(\tau)=\frac{1}{\ov{\tau}}+\mu, \ B(\tau)=\tau+2, 
$$
where $\hP_{y=0}$ is identified with $\hC$ via 
the map $\bsx=(x,0,z) \mapsto \tau=x+iz$. 
 (Here we use the notations $\ch G_\mu$ and $\ch A_\mu$ 
to distinguish them from 
$G_\mu$ and $A_\mu$ defined in 3.3.) 
Note that $\ch A_\mu$ is an orientation reversing 
conformal automorphism of $\hat \C$ 
and that $\ch A_\mu^2=A_\mu^2$. 
Let $\ch G_\mu^+$ be the index two  subgroup of $\ch G_\mu$
of orientation preserving transformations. 
Then we have 
\begin{eqnarray*}
\ch G^+_\mu=\langle \ch A_\mu^2, \,\ch A^{-1}_\mu B \ch A_\mu, \, B \rangle;  
\end{eqnarray*}
in fact, one can check that 
$\ch G_\mu^+ \subset \two$ and $[\ch G_\mu :\ch  G_\mu^+]=2$. 
Now observe that 
\begin{eqnarray*}
&&\ch A^{-1}_\mu B \ch A_\mu(\tau)=\frac{\tau}{2\tau+1}=C(\tau) \quad \text{and} \\
&&\ch A_\mu B \ch A^{-1}_\mu(\tau)=D_\mu(\tau)=C(\tau-\mu)+\mu.  
\end{eqnarray*}
It then follows 
that $H_\mu=\la B,C,D_\mu\ra$ is a subgroup of 
$\ch G_\mu^+=\la B,C, \ch A_\mu^2\ra$ and that $\ch G_\mu^+=\la H_\mu, \ch A_\mu^2\ra$. 
The statement of 
Theorem 5.4 can be rephrase that 
$\ch G_\mu^+$ is discrete provided that $H_\mu$ is discrete.  

To prove Theorem 5.4, we need to recall some basic facts of the action of a group 
$H_\mu=\la B,C,D_\mu \ra$ for $\mu \in \MN$ on $\hC$.  
For simplicity we assume that $\im \mu>0$, but 
the argument for the case $\im \mu <0$ is  parallel. 
Let $J_1=\la B,C \ra$ and $J_2=\la B,D_\mu \ra$ be the subgroups 
of $H_\mu$ of type $(0,3)$, and let 
$\Delta_1=\{\tau \in \C : \im \tau<0\}$ and 
$\Delta_2=\{\tau \in \C : \im \tau>\im \mu\}$. 
In this notation, we have the following: 

\begin{lem}[cf. \cite{Kr}]
\begin{enumerate}
  \item 
    Let $i=1,2$. 
    The disc $\Delta_i$ is a component of 
$\Omega(H_\mu)$, and is $(H_\mu,J_i)$-invariant;   
that is,  
  $h(\Delta_i)=\Delta_i$ for every $h \in J_i$ and 
  $h(\Delta_i) \cap \Delta_i=\emptyset$ 
for every $h \in H \setminus J_i$. 

\item $h(\Delta_1) \cap \Delta_2=\emptyset$ for every $h \in H_\mu$. 
\item 
The set 
$\Omega_0(H_\mu):=\Omega(H_\mu) \sm \bigcup_{h \in H_\mu} h(\Delta_1 \cup \Delta_2)$ 
is either an empty set or a connected component of $\Omega(H_\mu)$ which is 
$H_\mu$-invariant and simply connected. 
\end{enumerate}
\end{lem}

We state below the  second Klein--Maskit combination theorem  for Kleinian groups 
in $\three$ to what extent we need in the proof of Theorem 5.4 
(and its extension, Theorem 5.7).   
We refer the reader to Maskit \cite[VII]{Mas} for more information.  

\begin{thm}[The second Klein--Maskit combination theorem \cite{Mas}]
Let $H \subset \three$ be a torsion-free Kleinian group 
and $A \in \three$. 
Let $J_1,\,J_2$ be subgroups of $H$, and let 
$\B_1,\,\B_2 \subset \hR^3$ be closed topological balls. 
If they satisfy the following conditions $(1)$--$(5)$ then 
$G=\la H,A \ra$ is discrete and isomorphic to 
the HNN-extension $H *_A$ of $H$ by $A$: 
\begin{enumerate}
  \item   
  Let $i=1,2$. 
  The interior $\wa \B_i$ 
of $\B_i$ is $(H,J_i)$-invariant; that is,  
  $h(\wa \B_i)=\wa \B_i$ for every $h \in J_i$ and 
  $h(\wa \B_i) \cap \wa \B_i=\emptyset$ 
for every $h \in H \setminus J_i$. 
  \item 
  $h (\wa \B_1) \cap \wa \B_2=\emptyset$ for every $h \in H$. 
  \item The complement 
    $\hR^3 \setminus \bigcup_{h \in H}h(\B_1 \cup \B_2)$ 
  has an interior point. 
     \item The transformation $A$ takes 
     the interior of $\B_1$  
onto the exterior of $\B_2$; that is, 
$A(\wa \B_1) \cap \wa \B_2=\emptyset$ and 
  $A(\bd \B_1)=\bd \B_2$.  
  \item $J_2=AJ_1A^{-1}$. 
\end{enumerate}
\end{thm}

\begin{rem}
The conditions $(1)$, $(2)$ and $(3)$ guarantee that 
$\wa \B_1/J_1$ and $\wa \B_2/J_2$ can be embedded disjointly 
into $\Omega(H)/H$,  
and that the complement 
$\Omega(H)/H \sm (\wa \B_1/J_1 \cup \wa \B_2/J_2)$ 
has an interior point.  
The conditions $(4)$ and $(5)$ guarantee that 
the action of $A$ descends to 
a pairing map of the resulting boundary of 
$\Omega(H)/H \sm (\wa \B_1/J_1 \cup \wa \B_2/J_2)$. 
\end{rem}

For the convenience of the reader, we give here a sketch of the proof of Theorem 5.3; 
see \cite[VII. D and E]{Mas} for more details. 

\begin{proof}[Sketch of proof of Theorem 5.3] 

Recall from \cite{Mas} that the group $H*_A$ is the free group 
of words in $A$ and the elements of $H$ 
modulo equivalence induced from the relation $J_2=AJ_1A^{-1}$. 
Therefore 
each element of $H *_A$ is equivalent to a 
word of the form 
$A^{\alpha_n}h_n \cdots A^{\alpha_1}h_1$ 
($\alpha_i \in \Z$, $h_i \in H$) which satisfy the following conditions: 
(1) $h_i \ne \id$ for  $i>1$, (2) $\alpha_i \ne 0$ for $i<n$, 
 (3)  $\alpha_{i+1}<0$ if $\alpha_i<0$ and $h_{i+1} \in J_1 \sm \{\id \}$, 
and 
(4)  $\alpha_{i+1}>0$  if $\alpha_i>0$ and $h_{i+1} \in J_2 \sm \{\id \}$. 
We denote by $\Phi:H *_A \to \la H,A\ra \subset \three; \,
\Phi(A^{\alpha_n}h_n \cdots A^{\alpha_1}h_1)=A^{\alpha_n}h_n \cdots A^{\alpha_1}h_1$ 
the natural projection.

Now let us take an open set 
$U$ in $\hR^3 \sm \bigcup_{h \in H}h(\B_1 \cup \B_2)$. 
Since $H$ is a torsion-free Kleinian group, we may assume that 
$h(U) \cap U =\emptyset$ for every $h \in H \sm\{\id\}$. 
We will show that $\Phi(g)(U) \cap U=\emptyset$ for every 
$g=A^{\alpha_n}h_n \cdots A^{\alpha_1}h_1 \ne \id$. 
This implies that $\Phi$ is an isomorphism and $\la H,A \ra$ is discrete.

We first assume that $n=1$. 
If $\alpha_1=0$ then $g=h_1 \in H \sm \{\id\}$ and thus 
$g(U) \cap U=\emptyset$.  
When  $\alpha_1 \ne 0$, $g(U)=A^{\alpha_1}h_1(U)$ 
lies in $\wa \B_1$ if $\alpha_1<0$ and in $\wa \B_2$ if $\alpha_1>0$. 
In both cases we have $g(U) \cap U=\emptyset$.  

Next assume that $n=2$. 
If $\alpha_2=0$ then $g=h_2A^{\alpha_1}h_1$ and 
 $g(U) \subset \bigcup_{h \in H}h(\wa \B_1 \cup \wa \B_2)$, which implies that  
$g(U) \cap U =\emptyset$.  
When $\alpha_2 \ne 0$, we can show that 
$g(U)=A^{\alpha_2}h_2A^{\alpha_1}h_1(U) \subset \wa \B_1 \cup \wa \B_2$. 
In fact, if $h_2A^{\alpha_1}h_1(U) \subset \wa \B_1$, one see that 
$h_2 \in J_1$ and $\alpha_1<0$. 
This implies that $\alpha_2<0$ and that 
$g(U) \subset \wa \B_1$.  
Similarly, if $h_2A^{\alpha_1}h_1(U) \subset \wa \B_2$ we have 
$g(U) \subset \wa \B_2$. 
Finally, if $h_2A^{\alpha_1}h_1(U) \subset \hR^3 \sm (\wa \B_1 \cup \wa \B_2)$,   
$g(U)$ lies in $\wa \B_1$ if $\alpha_2<0$ and in $\wa \B_2$ if $\alpha_2>0$.  
Thus we have $g(U) \subset \wa \B_1 \cup \wa \B_2$, 
which implies that $g(U) \cap U =\emptyset$ also in this case. 
 
The proof for  the case of $n >2$ is obtained by induction. 
\end{proof}

Applying Lemma 5.2 and Theorem 5.3,  
we can now prove the following: 

\begin{thm}We have 
$$\hMM \cap P_{q=0}
=R_{\frac{\pi}{2}}({\mathcal M}_{0,4}).$$
\end{thm}

\begin{proof}
It follows from Theorem 5.1 that 
if $\bsp \in \hMM \cap P_{q=0}$ then 
 $\bsp \in R_{\frac{\pi}{2}}({\mathcal M}_{0,4})$. 
Conversely, suppose that  
$\bs p=(p,0,r) \in R_{\frac{\pi}{2}}(\MN)$, and hence 
that the subgroup $H_{\bs p}=\la B,C,D_\bsp\ra$ of 
$G_{\bs p}$ is discrete. 
We assume  for simplicity that the third coordinate $r$ of $\bsp=(p,0,r)$  
is positive, but the argument for the case $r<0$ is almost parallel. 
Let 
\begin{eqnarray*}
&&\B_1=\{(x,y,z) \in \R^3: z \le 0\} \cup \{\fty\}, \\
&&\B_2=\{(x,y,z) \in \R^3: z \ge r\} \cup \{\fty\}
\end{eqnarray*}
and 
\begin{eqnarray*}
J_1=\la B,C \ra, \quad 
J_2=\la B, D_\bsp \ra. 
\end{eqnarray*}
We will check that $H_\bsp,\,A_\bsp, \, \B_1,\,\B_2,\,J_1$ and $J_2$ 
satisfy the conditions $(1)$--$(5)$ in Theorem 5.3. 
Now take $\mu \in \MN$ such that $\bsp=R_\frac{\pi}{2}(\mu)$. 
Then  $\im \mu>0$ by our assumption $r>0$. 
The map $R_\frac{\pi}{2}:\hC \to R_\frac{\pi}{2}(\hC)=\hP_{y=0}$ 
conjugates the action of $H_\mu$ on $\hC$ to 
the action of $H_\bsp$ on the sphere $\hP_{y=0}$. 
Therefore the group $H_\bsp$ can be regarded 
as the Poincar\'{e} extension of the group $H_\mu$ acting on the 
sphere $\hP_{y=0} \cong \hC$ 
(see the left of Figure 6).  
Thus the conditions $(1)$ and $(2)$ 
directly follows from Lemma 5.2. 
Now observe that, for $i=1,\,2$, every  
orbit $h(\B_i)$  ($h \in H_\bsp$) of $\B_i$ except for $\B_i$ 
is a ball of radius $\le r/2$ with center in $P_{y=0}$. 
Therefore one can take an open domain, say,  
$\{(x,y,z) \in \R^3 : |x|>r/2, \, 0<z<r\}$ which does not intersect 
$\bigcup_{h \in H}h(\B_1 \cup \B_2)$. 
Thus the condition $(3)$ follows. 
It is easy to check the condition $(4)$. 
Finally, the condition $(5)$ $J_2=A_\bsp J_1A_\bsp^{-1}$ 
follows from the facts that  
$B=A_\bsp CA_\bsp^{-1}$ and $D_\bsp=A_\bsp BA_\bsp^{-1}$.  
Therefore Theorem 5.3 asserts that 
$G_\bsp=\la H_\bsp, A_\bsp \ra$ is discrete and isomorphic to the group 
$H_\bsp *_{A_\bsp}$.  
Thus we conclude that $\bsp \in \hMM$. 
\end{proof}

\begin{figure}[htbp]
  \begin{center}
\includegraphics[keepaspectratio=true,height=70mm]{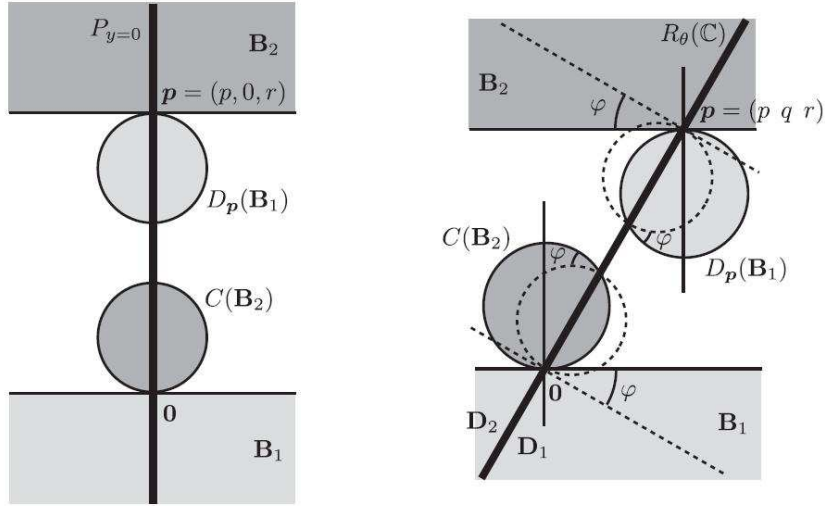}
\end{center}
  \caption{
  The balls $\B_1$, $\B_2$, $C(\B_2)$ and $D_\bsp(\B_1)$ 
  and the planes $P_{y=0}$ and $R_\theta(\C)$ in the proofs of 
  Theorems 5.2 (left) and 5.5 (right). 
  These are the views along the $x$-axis. }
\end{figure}

\subsection{Slices through the planes $R_\theta(\C)$}

Recall from Theorem 5.4 that 
the  Maskit slice $\MN$ appears as 
the slice of $\hMM$ through the plane 
$P_{q=0}=R_\frac{\pi}{2}(\hC)$. 
In this subsection, we will extend this result 
to Theorem 5.7, which states that there is a constant $\varphi_0>0$ such that 
$\MN$ also appears as the slice 
of $\hMM$ through the plane 
$R_\theta(\hC)$ for every 
$\theta \in [\pi/2-\varphi_0, \, \pi/2+\varphi_0]$. 

We first fix our terminology: 
\begin{defn}[lens] 
The intersection $K_1 \cap K_2$ 
of two spherical balls $K_1$, $K_2$ in $\hR^3$ 
is called a {\it lens}. 
The {\it inner angle} of lens $K_1 \cap K_2$ 
is the interior angle $\psi\,(0 <\psi<\pi)$ 
formed by two faces of $K_1 \cap K_2$. 
\end{defn}

The following lemma is  essential in the proof of Theorem 5.7; 
this is a special case of the collar lemma of Basmajian \cite[Theorem 1.1]{Ba}. 
 
\begin{lem}[Basmajian \cite{Ba}]
For a given $\mu \in \MN$, let  $\Delta$ 
be a non-invariant component of $\Omega(H_\mu)$. 
Let $\CC(\Delta)$ be the convex hull of $\Delta$ in $\HH^3$, 
and let $\Stab_{H_\mu}(\Delta)$  
denote the stabilizer of $\Delta$ in $H_\mu$. 
Then there is a constant $k_0>0$, 
which does not depend on the choices of $\mu \in \MN$ and 
$\Delta \subset \Omega(H_\mu)$,  
such that the $k_0$-neighborhood of $\CC(\Delta)$ in $\HH^3$ 
with respect to the hyperbolic metric is 
$(H_\mu,\Stab_{H_\mu}(\Delta))$-invariant. 
\end{lem}

\begin{rem}
Note that 
$\CC(\Delta) \subset \HH^3 \subset \hR^3$ 
can be regarded as a lens in $\hR^3$ with inner angle $\pi/2$.
Similarly, the $k_0$-neighborhood of $\CC(\Delta)$ 
in $\HH^3$ is a lens in $\hR^3$ 
with inner angle $\pi/2+\varphi_0$, where 
$\varphi_0$ is the constant determined by $k_0$.  
\end{rem}

We give here a proof of Lemma 5.6 for the convenience of the reader. 
\begin{proof}[Proof of Lemma 5.6]
Since $\Delta$ is non-invariant component of $\Omega(H_\mu)$, 
it follows form Lemma 5.2 that  
the stabilizer $\Stab_{H_\mu}(\Delta)$  of $\Delta$ 
is a group  of type $(0,3)$, and that 
$\Delta$ is $(H_\mu,\Stab_{H_\mu}(\Delta))$-invariant. 
Therefore we only need to show that 
there exists a constant $k_0>0$,  
which does not depend on the choices of $\mu$ and $\Delta$,  
such that the hyperbolic distance from 
$\CC(\Delta)$ to $h(\CC(\Delta))$ is greater than $2k_0$ 
for every $h \in H_\mu \sm \Stab_{H_\mu}(\Delta)$. 

Now let $\bd \CC(\Delta)$ denote the relative boundary of $\CC(\Delta)$ in $\HH^3$ and 
let $r:\HH^3 \to \bd \CC(\Delta)$ be the nearest point retraction. 
For a given $h \in H_\mu \sm \Stab_{H_\mu}(\Delta)$, let 
$U \subset \bd \CC(\Delta)$ denote the image 
of $h(\CC(\Delta))$ via the map $r$. 
Then one can check that $g(U) \cap U=\emptyset$ 
for every $g \in \Stab_{H_\mu}(\Delta)$. 
Therefore the hyperbolic area of $U$ is bounded above by 
the hyperbolic area ($=2\pi$) of the 
thrice-punctured sphere $\bd \CC(\Delta)/\Stab_{H_\mu}(\Delta)$. 
From this, we can deduce the existence of a desired constant $k_0$. 
\end{proof}

Using Lemma 5.6, we can now prove the following: 

\begin{thm}
There is a constant 
$0< \varphi_0 <\pi/2$ such that 
$$
\hMM \cap R_\theta(\C)=
R_\theta({\mathcal M}_{0,4})
$$
for every 
$\theta \in [\pi/2-\varphi_0, \, \pi/2+\varphi_0]$. 
\end{thm}

\begin{proof}
The argument is similar to the argument 
of the proof of Theorem 5.4. 
Let $0< \varphi_0 <\pi/2$ be the constant as in Remark of Lemma 5.6. 
We will show below that for every 
$\theta \in [\pi/2-\varphi_0, \, \pi/2+\varphi_0]$ we have 
$\hMM \cap R_\theta(\C)=R_\theta({\mathcal M}_{0,4})$. 
By symmetry, we may assume that $\pi/2 -\varphi_0 \le \theta \le \pi/2$. 

It follows from Theorem 5.1 that if  
 $\bsp \in \hMM  \cap R_\theta(\C)$ then 
 $\bsp \in R_\theta(\MN)$. 
Conversely, suppose that 
$\bsp=(p,q,r) \in R_\theta(\MN)$, 
and hence that the subgroup $H_\bsp=\la B,C,D_\bsp\ra$ of 
$G_\bsp$ is discrete. 
We again assume that $r>0$. 
Let $\B_1$, $\B_2 \subset \hR^3$ 
and $J_1, J_2 \subset H_\bsp$ be the 
same as in the proof of Theorem 5.4. 
We will check that $H_\bsp,\,A_\bsp, \, \B_1,\,\B_2,\,J_1$ and $J_2$ 
satisfy the conditions $(1)$--$(5)$ in Theorem 5.3. 
Since the conditions $(4)$ and $(5)$ 
are similarly satisfied,  
we only need to check the conditions $(1)$--$(3)$.  

Now take $\mu \in \MN$ such that $R_\theta(\mu)=\bsp$. 
Then $\im \mu>0$ by our assumption $r>0$. 
Observe that the action of group $H_\bsp$ 
preserves the sphere $R_\theta(\hat \C)$, 
and that the map $R_\theta:\hC \to R_\theta(\hC)$ conjugates 
the action of $H_\mu$ on $\hC$ to 
the action of $H_\bsp$ on the sphere $R_\theta(\hC)$. 
Therefore the group $H_\bsp$ can be 
regarded as the Poincar\'{e} extension of the group $H_\mu$ acting 
on the sphere $R_\theta(\hC) \cong \hC$ 
(see the right of Figure 6). 
Observe that $R_\theta$ takes 
the components $\Delta_1$, $\Delta_2$ of $\Omega(H_\mu)$ as in Lemma 5.2 
to the intersections 
\begin{eqnarray*}
\Delta_1':=\wa \B_1 \cap R_\theta(\hat \C) \quad \text{and} \quad 
\Delta_2':=\wa \B_2 \cap R_\theta(\hat \C)
\end{eqnarray*}
of the balls $\wa \B_1,\, \wa \B_2$ 
with the sphere $R_\theta(\hat \C)$, respectively. 
Therefore $\Delta_i'$ is $(H_\bsp, J_i)$-invariant for $i=1,2$, 
and $h(\Delta_1') \cap \Delta_2'=\emptyset$ for every $h \in H_\bsp$. 
Now let $\varphi$ denote $\pi/2-\theta$.   
Then our assumption 
$\pi/2 -\varphi_0 \le \theta \le \pi/2$ can be written as 
$0 \le \varphi \le \varphi_0$. 
Let $\D_1$, $\D_2$ be the two components of 
$\hR^3 \sm R_\theta(\hat \C)$ such that the inner angles of the lenses 
$\D_1 \cap \B_1$ and $\D_2 \cap \B_1$ are 
$\pi/2+\varphi$ and $\pi/2-\varphi$, respectively.  

(1) We first show that $\wa \B_1$ is $(H_\bsp,J_1)$-invariant. 
(The same argument reveals that $\B_2$ is 
$(H_\bsp,J_2)$-invariant.) 
It is easy to see that 
$h(\wa \B_1)=\wa \B_1$ for every $h \in J_1$. 
To show that  
$\wa \B_1 \cap h(\wa \B_1)=\emptyset$ for every $h \in H_\bsp \sm J_1$, 
it suffices to show that 
$(\D_1 \cap \wa \B_1) \cap (\D_1 \cap h(\wa \B_1))=\emptyset$. 
In fact,  $h$ takes the lens $\D_1 \cap \wa \B_1$ 
with inner angle  $\pi/2+\varphi$ 
 to the lens $\D_1 \cap h(\wa \B_1)$ with the same inner angle. 
 Since we are assuming that $0 \le \varphi \le \varphi_0$, 
 it follows from 
Lemma 5.6  that 
$(\D_1 \cap \wa \B_1) \cap (\D_1 \cap h(\wa \B_1))=\emptyset$. 

(2) We next show that
 $h(\wa \B_1) \cap \wa \B_2 =\emptyset$ for every $h \in H_\bsp$. 
 By symmetry, it suffices to show that 
 $(\D_1 \cap h(\B_1)) \cap (\D_1 \cap \B_2) =\emptyset$. 
 This follows from the facts that $h(\Delta_1') \cap \Delta_2' =\emptyset$, 
 that the lens
 $\D_1 \cap h(\B_1)$ has the inner angle $\pi/2+\varphi$,   
 and that the lens  $\D_1 \cap \B_2$ has the inner angle $\pi/2-\varphi$. 

(3) Finally we show that 
$\hR^3 \sm \bigcup_{h \in H_\bsp}h(\B_1 \cup \B_2)$ has an 
interior point. 
This  follows  from the fact that, for $i=1,\,2$, every  
orbit $h(\B_i)$  ($h \in H_\bsp$) of $\B_i$ except for $\B_i$ 
is a ball of radius $\le r/2$ which intersects with the plane 
$R_\theta(\hat \C)$. 

Therefore Theorem 5.3 asserts that 
$G_\bsp=\la H_\bsp, A_\bsp \ra$ is discrete and isomorphic to the group 
$H_\bsp *_{A_\bsp}$. 
Thus we conclude that $\bsp \in \hMM$. 
\end{proof}

\subsection{Slice thorough the plane $P_{p=0}$} 

In this subsection, we consider the case of $\bsp \in P_{p=0}$. 
In this  case the cyclic subgroup 
$\la A_\bsp \ra \subset G_\bsp$ preserves the sphere $P_{x=0}$ 
and its orientation. 
One of the main result of this subsection is Theorem 5.9,  
which gives a necessary and sufficient condition for $\bsp \in P_{p=0}$ 
to be contained in $\hMM$.  
As a consequence, we will show in Corollary 5.11 that the boundary of the slice of $\hMM$ 
through the plane 
$P_{p=0}$ is a union of countably many analytic arcs,  
although the boundaries  of the slices of $\hMM$ through the planes  
 $P_{r=0}$ and $P_{q=0}$ are so called ``fractal." 

We begin with the notion of a Ford domain: 
Let $\Gamma \subset \three$ be a Kleinian group and suppose that 
$\gamma(\fty) \ne \fty$ for every $\gamma \in \Gamma \sm \{\id\}$. 
Given $\gamma \in \Gamma$ we denote by 
$E(\gamma)$ the exterior of the isometric sphere 
$I(\gamma)$ of $\gamma$; that is, $E(\gamma)$ is 
the connected component of 
$\hR^3 \sm I(\gamma)$ containing $\fty$. 
Then the {\it Ford domain} 
for $\Gamma$ is defined by 
\begin{eqnarray*}
\mathrm{Ford}(\Gamma)=\bigcap_{\gamma \in \Gamma \sm \{\id\}}E(\gamma), 
\end{eqnarray*}
which turns out to be a fundamental domain for $\Gamma$. 
Below, we denote by $\radi(I(f))$ the radius of the isometric sphere 
$I(f)$ of a transformation $f \in \three$ with $f(\fty) \ne \fty$. 

We will need the following lemma in the proof of Theorem 5.9. 

\begin{lem}
Let $f \in \three$ be a loxodromic transformation with $f(\fty)\ne \fty$, 
and suppose that $f$ preserves the sphere $\hP_{x=0}$ 
and its orientation. 
Let $B(\bsx)=\bsx+(2,0,0)$. 
Then $[f,B]$ is loxodromic, pure parabolic  or elliptic
if and only if $\radi(I(f))<1,\,=1$ or $>1$, respectively. 
\end{lem}
 
\begin{proof}
Note that the points 
$f^{-1}(\fty)$, $f(\fty)$ lie in the plane $P_{x=0}$. 
We may assume that $f^{-1}(\fty)=\bs0$  after conjugating by a translation if necessary. 
We write $\bsv=f(\fty)$.  
Then one see from 2.3 that  
$f(\bsx)=PJ_{I(f)}(\bsx)+\bsv$ for some 
$P \in O(3) \sm SO(3)$. 
Moreover, since $f$ preserves the sphere $\hP_{x=0}$ and its orientation,  
one can see that $P(\bse_1)=\bse_1$, where $\bse_1=(1,0,0)$. 

Now write $r=\radi(I(f))$ and suppose first that $r=1$. 
Then we see that 
the transformation $[f,B]=fBf^{-1}B^{-1}$ maps 
the  exterior of the sphere 
$S_1:=B(I(f^{-1}))$ onto 
the interior of the sphere $S_2$ of radius $1/4$ with center at 
$\bsv+3/4\,\bse_1$ (see Figure 7).  
Note that $S_1$ touches to $S_2$ at $\bsv+\bse_1$, and that 
the transformation $[f,B]$ fixes the point $\bsv+\bse_1$. 
Furthermore, we see that the differential 
of the map $[f,B]$ at this point is the identity. 
Thus we conclude that $[f,B]$ is pure parabolic.  
By the same argument, we see that 
if $r<1$ then $S_1 \cap S_2=\emptyset$ and thus $[f,B]$ is loxodromic, 
and if $r>1$ then $[f,B]$ fixes every points of $S_1 \cap S_2$ 
and thus $[f,B]$ is elliptic. 
\end{proof}

\begin{figure}[htbp]
\begin{center}
\includegraphics[keepaspectratio=true,height=60mm]{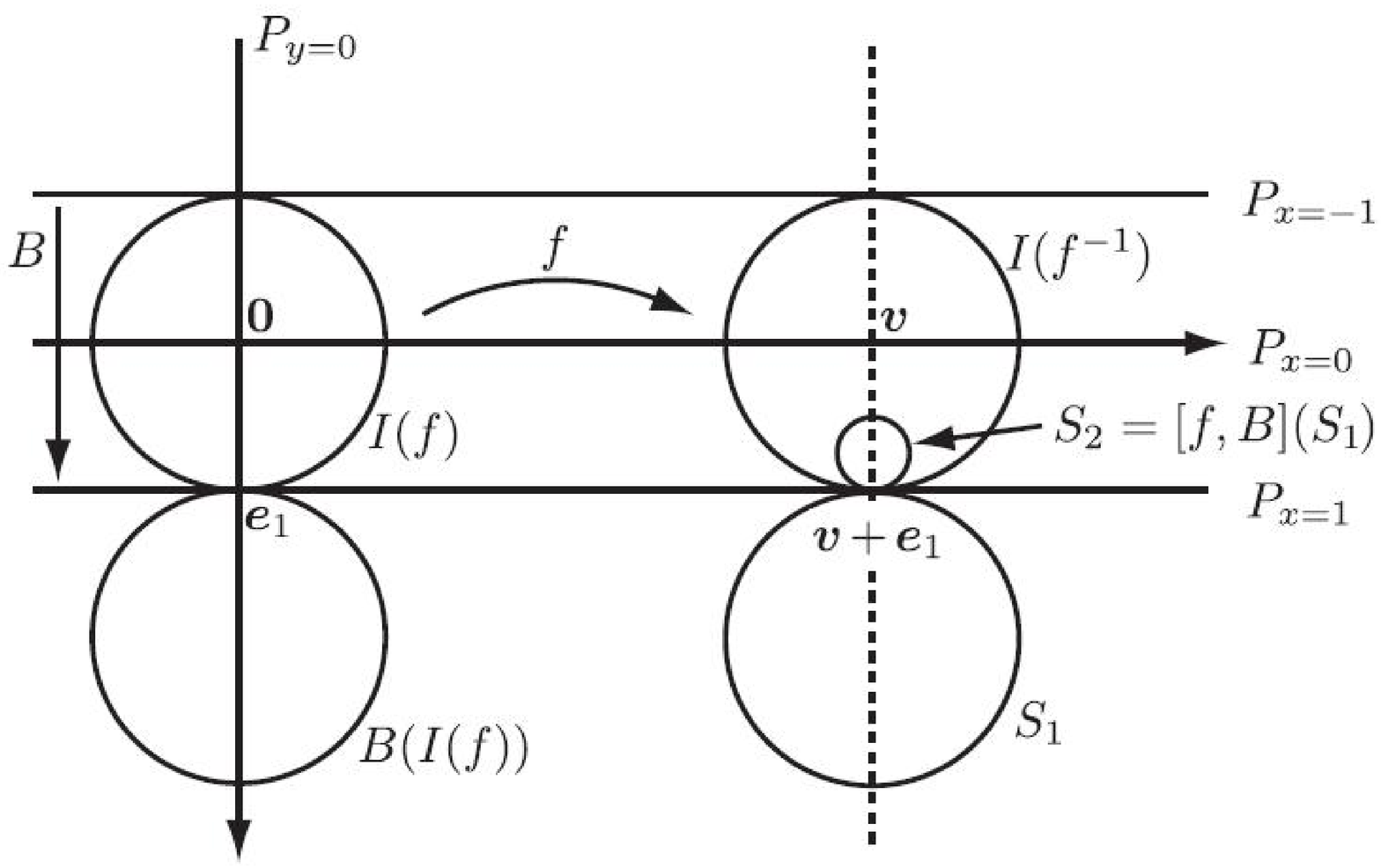}
\end{center}
  \caption{Surfaces $S_1$ and $S_2$ in the case of $r=1$. 
  This is the view along the $z$-axis.}
\end{figure}

\begin{thm}
Let $\bsp \in P_{p=0}$. 
Then $\bsp \in \hMM$ if and only if 
$\radi(I(A_\bsp^n)) \le 1$ for all $n \in \N$.   
Moreover, for a given $n \in \N$, $\radi(I(A_\bsp^n))=1$ 
if and only if $[A_\bsp^n,B]$ is pure parabolic. 
\end{thm}

\begin{proof}
Let $\bsp \in P_{p=0}$. 
We first remark that $A_\bsp^n(\fty) \ne \fty$ 
for every $n \in \Z\sm\{0\}$, because 
$A_\bsp^n(\fty)$ lies in 
the interiors of $I(A_\bsp)$ or $I(A_\bsp^{-1})$.  
Therefore the Ford domain 
$$
\mathrm{Ford}(\la A_\bsp\ra)=\bigcap_{n \in \Z \sm \{0\}}E(A_\bsp^n)
$$
for the cyclic group $\la A_\bsp \ra$ can be defined. 
In addition, observe that 
the center $A_\bsp^n(\fty)$ of $I(A_\bsp^n)$ 
lies in the plane $P_{x=0}$ for every $n \in \Z \sm \{0\}$,  
and by definition that 
$\radi(I(A_\bsp^n))=\radi(I(A_\bsp^{-n}))$ for every $n \in \N$. 

Suppose first that $\radi(I(A_\bsp^n)) \le 1$ for all $n \in \N$. 
Then the Ford domain $\mathrm{Ford}(\la A_\bsp\ra)$ 
contains the set $\{(x,y,z) \in \R^3:|x| \ge 1\}$. 
Therefore we can apply the Klein--Maskit combination theorem 
(Theorem 5.3) to show that $G_\bsp=\la A_\bsp, B \ra$ is discrete 
and isomorphic to the rank-$2$ free group $\la A_\bsp \ra *_B$ 
by letting $\B_1=\{(x,y,z) \in \R^3 : x \le -1\}$, 
$\B_2=\{(x,y,z) \in \R^3 : x \ge 1\}$ and $J_1=J_2=\{\id\}$. 

Conversely, suppose that  $\radi(I(A_\bsp^n))>1$ for some $n$. 
Then by Lemma 5.8, $[A_\bsp^n,B]$ is elliptic. 
If the order of this elliptic element is infinite then the group $G_\bsp$ is non-discrete, 
and if the order is finite then $G_\bsp$ is not free. 
In both cases, we conclude $\bsp \not\in \hMM$. 

The second  statement also follows from Lemma 5.8. 
\end{proof}
Now let $\bsp=(0,q,r) \in P_{p=0}$ and $\mu=q+ir \in \C$. 
In this notation, 
the map $A_\bsp(\bsx)=\hat JJ(\bsx)+\bsp$ 
restricted to the sphere $\hP_{x=0}$ is a M\"{o}bius transformation given by 
$$
\ac A_\mu(\tau)=-\frac{1}{\tau}+\mu, 
$$
where $\hP_{x=0}$ is identified with $\hC$ via the map 
$\bsx=(0,y,z) \mapsto  \tau =y+iz$. 
(Here we use the notation $\ac A_\mu$ 
to distinguish it from $A_\mu$ and $\ch A_\mu$.) 
The matrix representation of $\ac A_\mu$ is 
\begin{eqnarray*}
\ac A_\mu=\left(
  \begin{array}{cc}
    \mu   &   -1 \\
    1   &  0  \\
  \end{array}
\right) \in \psl. 
\end{eqnarray*}
We denote the entries of $\ac A_\mu^n$ for $n \in \N$ by 
\begin{eqnarray*}
\ac A_\mu^n=\left(
  \begin{array}{cc}
    a_n   &b_n    \\
     c_n  &d_n    \\
  \end{array}
\right). 
\end{eqnarray*}
Then the radius  
of the isometric sphere $I(\ac A_\mu^n)$  of $\ac A_\mu^n$ is given by $1/|c_n|$.   
In other words, we have 
$\radi(I(A_\bsp^n))=\radi(I(A_\bsp^{-n}))=1/|c_n|$ 
for every $n \in \N$.  
Since 
\begin{eqnarray*}
\left(
  \begin{array}{cc}
    a_{n+1}   & b_{n+1}   \\
    c_{n+1}   &  d_{n+1}  \\
  \end{array}
\right)&=&
\left(
  \begin{array}{cc}
    \mu   &   -1 \\
    1   &  0  \\
  \end{array}
\right)
\left(
  \begin{array}{cc}
    a_n   &b_n    \\
     c_n  &d_n    \\
  \end{array}
\right) \\
&=&\left(
  \begin{array}{cc}
    \mu a_n-c_n   &   \mu b_n-d_n \\
    a_n   &  b_n \\
  \end{array}
\right), 
\end{eqnarray*}
we have 
$a_{n+1}=\mu a_n-c_n$ and $c_{n+1}=a_n$. 
Therefore we obtain the following recurrence equations: 
$$
c_1=1, \quad  c_2=\mu \quad \text{and} \quad c_{n+2}=\mu c_{n+1}-c_n \ (n\in \N). 
$$
Especially $c_n=c_n(\mu)$ is a monic $\mu$-polynomial of degree $(n-1)$; 
for example  
\begin{eqnarray*}
&&c_3(\mu)=\mu^2-1, \\
&&c_4(\mu)=\mu^3-2\mu, \\
&&c_5(\mu)=\mu^4-3\mu^2+1, 
\end{eqnarray*}
and so on. 
Using this notation and identifying 
$P_{p=0}$ with $\C$, we can rephrase Theorem 5.9 as follows.  

\begin{thm}
We have 
$$
\hMM \cap P_{p=0}=\bigcap_{n \in \N}\{\mu \in \C: 
|c_n(\mu)| \ge1\}. 
$$
\end{thm}

Since the set $\{\mu \in \C : |c_n(\mu)|=1\}$ 
is a  one-dimensional real analytic variety 
for every $n \in \N$, we have the following corollary.  

\begin{cor}
The boundary of $\hMM \cap P_{p=0}$ in the plane $P_{p=0}$ 
is a union of countably many analytic arcs. 
\end{cor}

One can find in Figures 8 and 9 the numerical graphics of the 
loci of $\mu \in \C$ such that $|c_n(\mu)|=1$ or $ \le 1$ for some $n$. 
\begin{figure}[htbp]
  \begin{center}
\includegraphics[keepaspectratio=true,height=45mm]{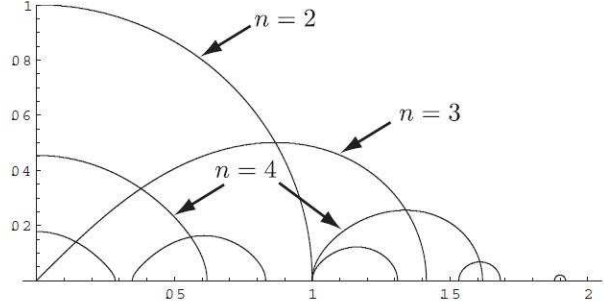}
  \end{center}
  \caption{The loci of 
  $\mu \in \{q+ir \in \C : 0 \le q \le 2, \, 0 \le r \le 1\}$ 
   such that 
  $|c_n(\mu)|=1$ for $n=2,3,4$ and $10$.  
 The arcs without labels correspond to $n=10$.}
\end{figure}

\begin{figure}[htbp]
  \begin{center}
    \includegraphics[keepaspectratio=true,height=45mm]{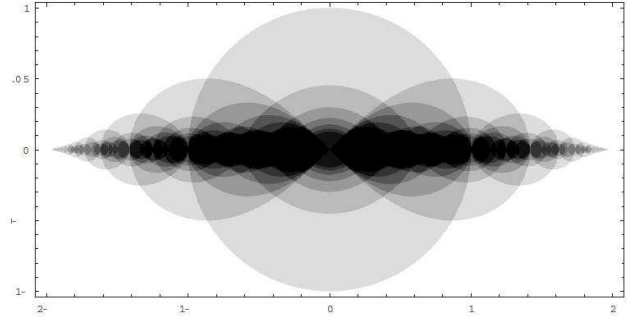}
  \end{center}
  \caption{The union of the loci of 
  $\mu \in \{q+ir \in \C : -2 \le q \le 2, \, -1 \le r \le 1\}$ 
  such that 
  $|c_n(\mu)| \le 1$ for $2 \le n \le 20$.   
  The darkness increases as the number of the intersection increases.}
\end{figure}

\section{Appendix}

In this section we consider a family of $3$-generator Kleinian groups in $\three$, 
which can be viewed as an analogue of the family of groups of type $(1,1)$ in $\two$. 
These groups contain groups of type $(1,1)$ in $\three$ and their 
limit sets are union of round spheres. 

For a given $\bsp=(0,q,0) \in \R^3$ with $q > 2$, 
we define an ideal hexahedron $\DD_\bsp$  as follows (see Figure 10): 
\begin{eqnarray*}
 \DD_\bsp=
 \left\{
\bsx=(x,y,z)\in \R^3 : 
  \begin{aligned}
   &|x| \le 1/\sqrt{2}, \, 0 \le y \le q, \, |z| \le 1/\sqrt{2},     \\
   & |\bsx| \ge 1,  \, |\bsx-\bsp| \ge 1   \\
  \end{aligned}  
\right\}. 
\end{eqnarray*}
Observe that eight edges of $\DD_\bsp$ have dihedral angle $\pi/4$ and 
the rest four edges have dihedral angle $\pi/2$. 
Moreover, the following three transformations in $\three$ pair the faces 
of $\DD_\bsp$: 
\begin{eqnarray} \quad \quad \quad 
A_\bsp(\bsx)= \hat J J(\bsx)+\bsp, \  
B(\bsx)=\bsx+(\sqrt{2},0,0),  \ 
C(\bsx)=\bsx+(0,0,\sqrt{2}). 
\end{eqnarray}
(Although the notation $A_\bsp$ is common with the previous sections, 
we remark that $B$, $C$ are different from those in the previous sections.)  
Therefore, it follows from 
Poincar\'{e}'s polyhedron theorem (cf. \cite{Mas}) that the group 
$$
K_\bsp=\la A_\bsp, B,C \ra
$$ 
generated by $A_\bsp,\,B$ and $C$ in $\three$
is discrete, 
and that $K_\bsp$ is isomorphic to the  abstract group 
$$
\mathfrak{g}=\la a,b,c:[a,b]^2=[a,c]^2=[b,c]=\id \ra. 
$$ 

For arbitrary $\bsp \in \R^3$, we also denote by 
$K_\bsp=\la A_\bsp,B,C \ra$ the group in $\three$ generated by the 
transformations $A_\bsp$, $B$ and $C$ as in $(6.1)$. 
Note that if $K_\bsp$ is isomorphic to $\mathfrak{g}$, 
the subgroup $\la A_\bsp, BC \ra$ of $K_\bsp$ 
is a group of type $(1,1)$ in $\three$ 
because $BC(\bsx)=\bsx+(\sqrt{2},0,\sqrt{2})$ is a translation of 
length $2$. 

\begin{figure}[t]
  \begin{center}
    \includegraphics[keepaspectratio=true,height=50mm]{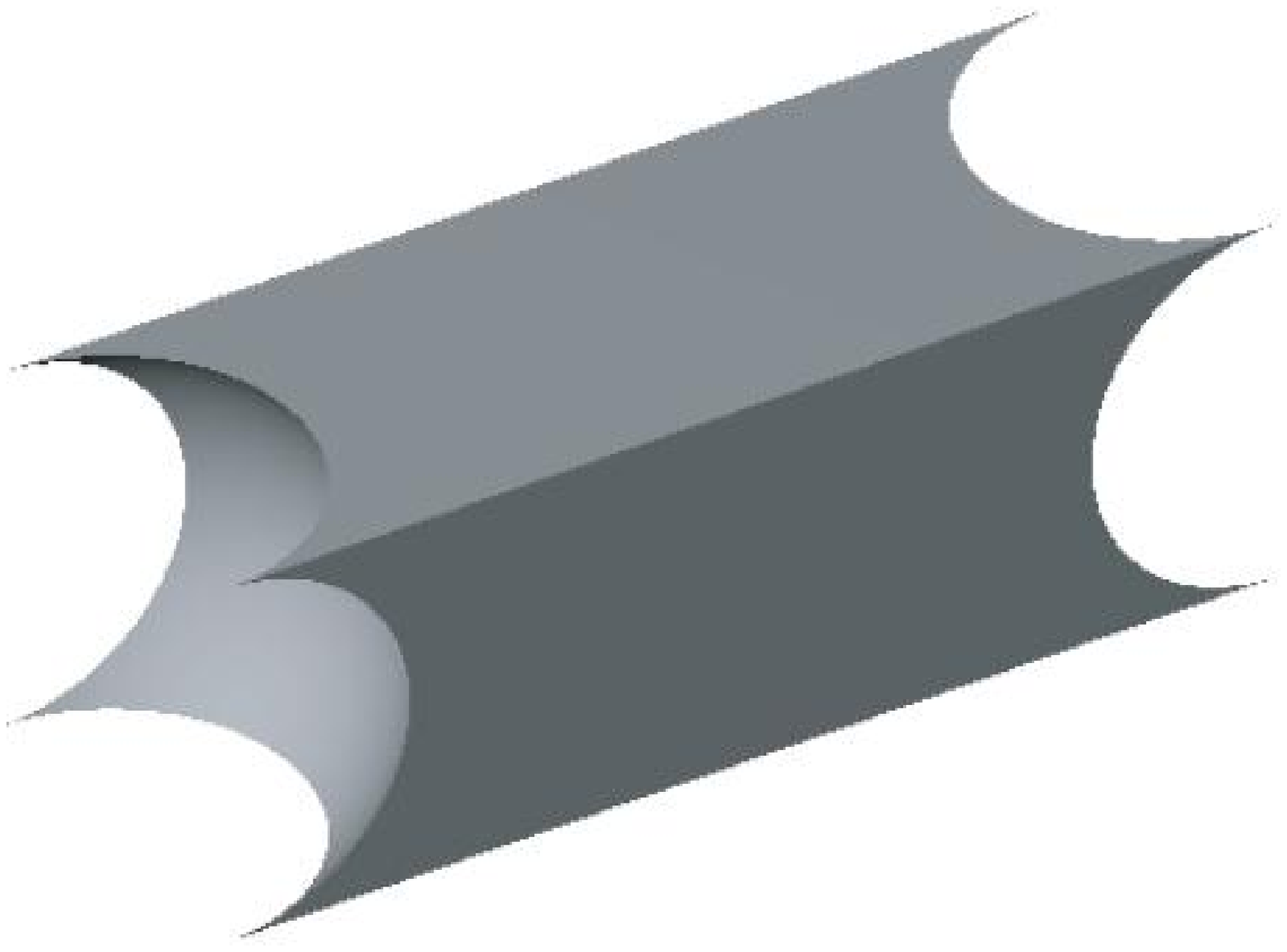}
  \end{center}
  \caption{The ideal hexahedron $\DD_\bsp$ with $\bsp=(0,4,0)$. }
\end{figure}

\subsection{Groups of type $(1,0;2)$} 
For a given $\bsp_0=(0,q,0) \in \R^3$ with $q > 2$,  
we will study deformations of $K_{\bsp_0}$ in the next subsection. 
As a preparation for this purpose, 
in this subsection, we study deformations of the subgroup 
$\la A_{\bsp_0}, B \ra$ of $K_{\bsp_0}=\la A_{\bsp_0},B,C \ra$.  

Let $\Si_{1,0;2}$ denote a torus with a singular point of cone-angle $\pi$. 
The fundamental group $\pi_1(\Si_{1,0;2})$ of the orbifold 
$\Si_{1,0;2}$ is isomorphic to an abstract group 
$\la a,b:[a,b]^2=\id \ra$. 
We say that a group $\Gamma$ in $\three$ is of {\it type} $(1,0;2)$ 
if it is the image of a faithful representation 
$\rho:\pi_1(\Si_{1,0;2})=\la a,b:[a,b]^2=\id \ra \to \three$ 
such that $\rho(b)$ is pure parabolic.  
Note that if $K_\bsp=\la A_\bsp, B, C \ra$ is isomorphic to $\mathfrak{g}$, 
the subgroups $\la A_\bsp, B \ra$, $\la A_\bsp, C \ra$  
of $K_\bsp$ are of type 
$(1,0;2)$ in $\three$. 
Similar to Theorem 4.5, we have the following normalization of 
groups of type $(1,0;2)$ in $\three$. 

\begin{thm}
Let $\Gamma=\la \alpha, \beta \ra$ be a subgroup of $\three$  
which is isomorphic to the abstract group 
$\la a,b:[a,b]^2=\id \ra$ and suppose that 
$\beta$ is pure parabolic.  
Then $\Gamma$ is conjugate in $\three$ to $\la A_\bsp,B \ra$ 
for some $\bsp=(p,q,r) \in \R^3$, 
where $A_\bsp, B \in \three$ are as in $(6.1)$. 
\end{thm}

\begin{proof}
The proof is almost parallel to the proof of Theorem 4.5. 
We can  first show that $\Fix(\alpha) \cap \Fix(\beta)=\emptyset$. 
In fact, if not, the same argument as in the proof of Theorem 4.5 implies that 
$[\beta,[\alpha,\beta]]=\id$, which 
contradicts to the assumption that 
$\Gamma=\la \alpha, \beta \ra$ is isomorphic to the abstract group 
$\la a,b:[a,b]^2=\id \ra$. 
Therefore,  by using the same normalization as in the proof of Theorem 4.5, 
we may assume that $\alpha,\beta$ are of the forms 
$\alpha(\bsx)=PJ(\bsx)+\bsp$, $\beta(\bsx)=\bsx+(u,v,0)$ 
where $P \in O(3) \sm SO(3)$, $\bsp \in \R^3$ and $u,v \ge 0$ 
such that $P^{-1}(\bsu)=(u,-v,0)$. 
Then $\alpha^{-1}\beta^{-1}\alpha\beta$ preserves the sphere $\hP_{z=0} \cong \hC$ and 
its orientation. 
Furthermore, identifying $\hP_{z=0}$ with $\hC$ via the map $(x,y,0) \mapsto x+iy$ and 
letting $\mu=u+iv$, one see that 
the map 
$\alpha^{-1}\beta^{-1}\alpha\beta$ restricted to the sphere $\hP_{z=0} \cong \hC$ 
is a M\"{o}bius transformation 
whose matrix representation is 
\begin{eqnarray*}
\left(
  \begin{array}{cc}
     1  & \mu   \\
     -\mu  &  1-\mu^2  \\
  \end{array}
\right) \in \psl. 
\end{eqnarray*}
It then follows from the condition $[\alpha,\beta]^2=\id$ 
that $\mu=\sqrt{2}$, and hence, that $\bsu=P^{-1}(\bsu)=(\sqrt{2},0,0)$. 
The remaining argument is the same to the argument of the proof of Theorem 4.5. 
\end{proof}

\subsection{Deformations of the group $K_\bsp$}

For a given $\bsp_0=(0,q,0) \in \R^3$ with $q > 2$, 
we now consider deformations $\{\phi:K_{\bsp_0} \to 
\Gamma \}$ of $K_{\bsp_0}$ in $\three$ 
such that the isomorphism 
$\phi$ takes a pure parabolic transformation to a pure parabolic transformation. 
The following theorem reveals that every such deformation $\Gamma$ of $K_{\bsp_0}$ 
is conjugate to $K_\bsp$ for some $\bsp \in \R^3$. 

\begin{thm}
Let $\Gamma=\la \alpha, \beta,\gamma \ra$ be a subgroup of $\three$ 
which is isomorphic to the abstract group  
$\mathfrak{g}=\la a,b,c:[a,b]^2=[a,c]^2=[b,c]=\id \ra$, 
and suppose  that $\beta,\gamma$ and $[\alpha,\beta\gamma]$ are pure parabolic. 
Then $\Gamma$ is conjugate in $\three$ to 
$K_\bsp=\la A_\bsp, B, C \ra$ for some $\bsp \in \R^3$. 
\end{thm}

\begin{proof}
Note that the subgroups $\la \alpha, \beta \ra$, $\la \alpha, \gamma \ra$ of $\Gamma$ 
are of type $(1,0;2)$ in $\three$, 
and $\la \alpha, \beta\gamma \ra$ is of type $(1,1)$.  
By Theorem 6.1, we may assume that 
$\la \alpha,\beta \ra=\la  A_\bsp,B \ra$ for some $\bsp \in \R^3$ 
after conjugating $\la \alpha,\beta,\gamma \ra$ if necessary.  
We will show that $\gamma=C^{\pm 1}$ below.  
Since  $\gamma$ is pure parabolic and commutes with $B$, 
we have $\gamma(\infty)=\infty$. 
Since we are normalizing  so that the radius of the isometric sphere of 
$\alpha=A_\bsp$ equals $1$, it follows from 
Theorem 6.1 that $\gamma$ is a translation of length $\sqrt{2}$,  
and from Theorem 4.5 that $B\gamma$ is a translation of length $2$. 
Therefore the direction of translations of $B$ and $\gamma$ 
are perpendicular. This implies that $\gamma=C^{\pm1}$. 
Thus we obtain $\la \alpha,\beta,\gamma \ra=\la A_\bsp, B,C \ra$. 
\end{proof}

We now define:
$$
\NN=\{\bsp=(p,q,r) \in \R^3 : K_\bsp=\la A_\bsp,B,C \ra \text{ is discrete and } \cong \mathfrak{g} \}. 
$$
Since the group $K_\bsp$ for $\bsp \in \NN$ contains 
the groups $\la A_\bsp, B^{-1}C \ra$, $\la A_\bsp, BC \ra$ of type $(1,1)$ in $\three$, 
the set $\NN\subset \R^3$ lies in the rotations of the set 
$\hMM$ of angles $\pi/4,\,3\pi/4$ along 
the $q$-axis.  
One can also see that $\NN$ is invariant under the action of the translations 
$\bsx \mapsto \bsx+(\sqrt{2}\,m,0, \sqrt{2}\,n)$, $m,n \in \Z$. 
The following lemma shows that the set $\NN$ contains a $3$-dimensional domain of $\R^3$. 

\begin{lem} We have 
$\{\bsp=(p,q,r) \in \R^3: |q| \ge 2\} \subset \NN$. 
\end{lem}

\begin{proof}
For a given $\bsp=(p,q,r) \in \R^3$ with $|q| \ge 2$, we let $\bsp_0=(0,q,0)$. 
We suppose that $q \ge 2$ for simplicity, but 
the argument for the case of $q  \le -2$ is almost parallel. 
Deforming the ideal hexahedron $\DD_{\bsp_0}$ in a similar way as in Figure 4,  
we can obtain an ideal polyhedron 
which can be applied Poincar\'{e}'s polyhedron theorem (cf. \cite{Mas}) 
to show  the group $K_\bsp=\la A_\bsp,B,C \ra$ 
is discrete and isomorphic to 
the abstract group $\mathfrak{g}$. 
Thus we conclude that $\bsp \in \NN$. 
\end{proof}

\subsection{Limit sets of Kleinian groups $K_\bsp$}

For a given $\bsp \in \R^3$, one can see that 
$A_\bsp^{-1}BA_\bsp=JBJ$ and $A_\bsp^{-1}CA_\bsp=JCJ$ hold. 
Therefore the group $K_\bsp=\la A_\bsp, B,C \ra$ always contains the group 
$$
L:=\la B,C,JBJ,JCJ \ra.  
$$
Note that the action of $L$ preserves the sphere 
$\hP_{y=0}$ and its orientation.  
Now we take a ideal pentahedron $\DD^-$ as follows:  
\begin{eqnarray*}
\DD^-=
\{\bsx=(x,y,z)\in \R^3 : 
|x| \le 1/\sqrt{2}, \, |z| \le 1/\sqrt{2}, \, 
|\bsx| \ge 1, \,  y \le 0 \}.  
\end{eqnarray*} 
Then one can see that the ideal octahedron 
$\DD^- \cup J(\DD^-)$ 
is a fundamental domain for the action of $L$ acting on the half-space 
$\{(x,y,z) \in \R^3 : y <0\}$. 
Therefore $L$ can be regarded as a Kleinian group in $\two$ of the first-kind, 
that is, the limit set $\Lambda(L)$ of $L$ 
is the whole sphere $\hP_{y=0}$. 
Thus the limit set $\Lambda(K_\bsp)$  of a Kleinian group $K_\bsp$ for $\bsp \in \NN$ 
is the closure of a union of round spheres: 
\begin{eqnarray*}
\Lambda(K_\bsp)=\overline{\bigcup_{\gamma \in K_\bsp/L}\gamma(\Lambda(L))}. 
\end{eqnarray*}
One can find in Figure 11 some computer graphics of the limit sets of $K_\bsp$ 
for parameters $\bsp \in \R^3$ which lie (or seem to be lie) in $\NN$.

By using the observation of the subgroup $L \subset K_\bsp$ above,   
we obtain an alternative proof of Lemma 6.3. 
In fact, a Kleinian group $K_\bsp$ for $\bsp=(p,q,r)$ with $|q| \ge 2$ 
can be regarded  as an amalgamation 
of the Kleinian group $L$ with the cyclic group 
$\la A_\bsp \ra$ as follows: 
Suppose $q \ge 2$ for simplicity. 
Then we can apply the second Klein--Maskit combination theorem (Theorem 5.3)
to show that $K_\bsp=\la L,A_\bsp \ra$ is discrete 
and $K_\bsp \cong L*_{A_\bsp}$ by letting 
\begin{eqnarray*}
&&\B_1=\left\{\bsx \in \R^3 : \left|\bsx-\left(0,\frac{q}{4},0\right)\right| 
\le \frac{q}{4}\right\}, \\
&&\B_2=\left\{(x,y,z) \in \R^3 : y \ge q-\frac{2}{q} \right\}, 
\end{eqnarray*}
$J_1=\la JBJ,JCJ \ra$ and $J_2=\la B,C \ra$.

\begin{figure}[htbp]
  \begin{center}
    \includegraphics[keepaspectratio=true,height=120mm]{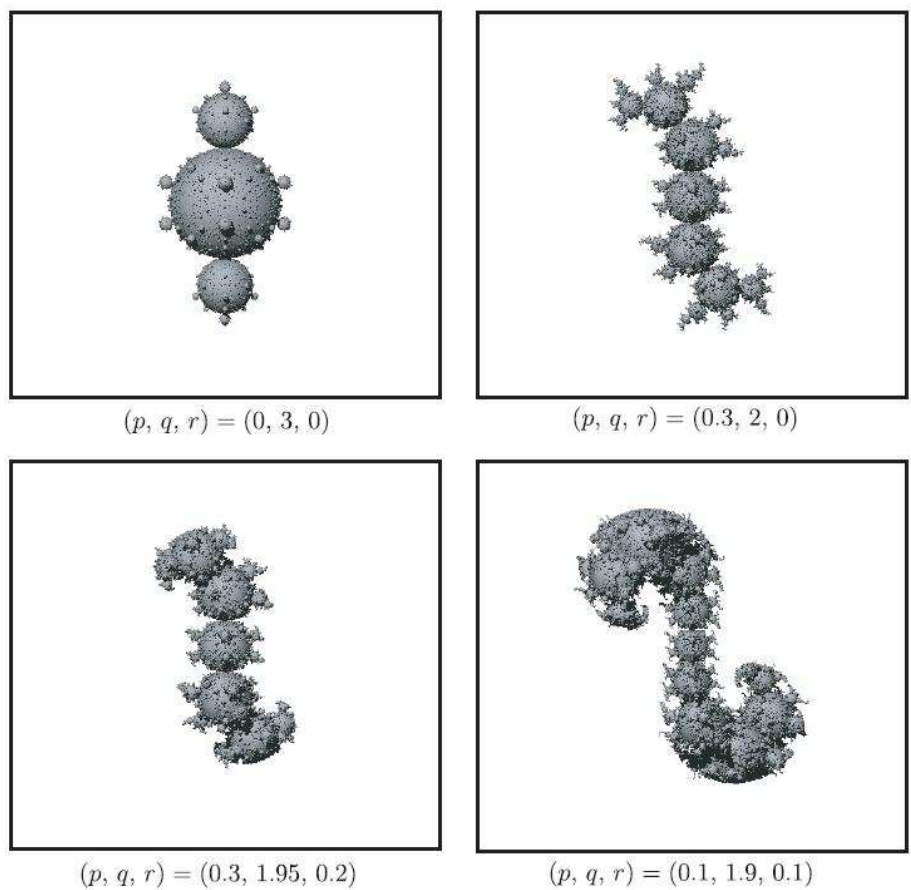}
  \end{center}
  \caption{Computer graphics of the limit sets $\Lambda(K_\bsp)$ 
of groups $K_\bsp$:  
These are the images $f(\Lambda(K_\bsp))$ of $\Lambda(K_\bsp)$ 
by the M\"{o}bius transformation 
$f(\bsx)=2 \hat J J(\bsx-\bse_2)-\bse_2$, $\bse_2=(0,1,0)$,  
which takes the half-space $y \le 0$ to the unit ball $|\bsx| \le 1$ and 
the point $\bse_2$ to the infinity. 
These  are the views along the $z$-axis.} 
\end{figure}

\end{document}